\documentclass[12pt]{article}
\usepackage{latexsym}
\usepackage{amsfonts}
\def\Bbb{\mathbb}
\def\bcp{\mathbb C \mathbb P}
\def\eea{\end{eqnarray*}}
\def\mv{\mbox{Vol}}

\newtheorem{defn}{Definition}
\newtheorem{thm}{Theorem}[section]
\newtheorem{prop}[thm]{Proposition}
\newtheorem{cor}[thm]{Corollary}
\newtheorem{lem}[thm]{Lemma}
\newtheorem{conj}[thm]{Conjecture}
\newtheorem{question}[thm]{Question}
\newenvironment{proof}{\medskip \noindent
{\bf Proof.}}{\hfill \rule{.5em}{1em}
\\}
\newenvironment{xpl}{\mbox{ }\\ {\bf  Example}\mbox{ }}{
\hfill $\diamondsuit$\mbox{}\bigskip}
\newenvironment{rmk}{\mbox{ }\\{\bf  Remark}\mbox{ }}{
\hfill $\Box$\mbox{}\bigskip}
\begin{document}
\sloppy
\title{Ricci Curvature, Minimal Volumes,\\ and Seiberg-Witten Theory}

\author{Claude LeBrun\thanks{Supported 
in part by  NSF grant DMS-9802722.} 
\\ 
SUNY Stony
 Brook 
  }

\date{March 9, 2000}
\maketitle

\begin{abstract} We derive new, sharp lower bounds for certain curvature
functionals on the space of Riemannian metrics of  
 a smooth compact $4$-manifold with a non-trivial
Seiberg-Witten invariant. These allow one, for example, 
to exactly compute the infimum of the $L^{2}$-norm of
 Ricci curvature for all  complex surfaces of general type. 
We are also   able to  show that the standard metric
on any 
complex hyperbolic $4$-manifold  minimizes volume among 
all metrics satisfying a point-wise lower bound on  sectional curvature
plus suitable multiples of the scalar curvature. 
These estimates also imply new non-existence results for
Einstein metrics. 
 \end{abstract}
%\vfill
%\pagebreak

\section{Introduction}

The aim of modern Riemannian geometry 
is to understand  the relationship between  
topology and curvature. 
It is thus  interesting and natural to 
consider differential-topological invariants of a smooth 
compact $n$-manifold which, by their very definition,
represent quantitative obstructions to the 
existence of  a  scalar-flat
 (or Ricci-flat, or flat)  
 metric on the given manifold. 
For example, if $M$ is a smooth compact $n$-manifold,
one may define invariants 
\begin{eqnarray*}
	{\mathcal I}_{s}(M) & = & 
 	\inf_{g}\int_{M}|{s}_{g}|^{n/2}d\mu_{g} \\
	{\mathcal I}_{r}(M) & = & 
	\inf_{g}\int_{M}|{r}|^{n/2}_{g}d\mu_{g}  \\
	{\mathcal I}_{\mathcal R}(M) & = & 
	\inf_{g}\int_{M}|{\mathcal R}|^{n/2}_{g}d\mu_{g} 
\end{eqnarray*}
where the infima are to be taken over all metrics
$g$ on $M$,  where  $s$, $r$, and  $\mathcal R$ 
respectively denote the 
scalar, Ricci, and  Riemann 
curvatures of $g$, and where $d\mu$ denotes 
Riemannian volume measure.
The power of $n/2$ is used for reason of  scale invariance: 
any other choice would result in the zero invariant.  
Still, one has every right to suspect that 
invariants with such  soft definitions
 will  either be trivial, or else 
will  be completely impossible to calculate
in practice. 

This suspicion would seem to be vindicated 
by a recent result of Petean \cite{jp3}, who, building upon 
the earlier work of Gromov-Lawson \cite{gvln} and Stolz
\cite{stolz},
showed that ${\mathcal I}_{s}(M)=0$ for any 
simply connected 
$n$-manifold, $n \geq 5$.
Dimension $4$, however, turns out to be radically different. 
Seiberg-Witten theory naturally leads to 
non-trivial lower bounds for ${\mathcal I}_{s}$
which,  amazingly, are often sharp. Using this,
the author has elsewhere \cite{lno,lky} computed ${\mathcal I}_{s}(M)$
for all complex surfaces $M$ with even first Betti number; 
it turns out that ${\mathcal I}_{s}(M)$  is positive
exactly for the surfaces of general type, and for
these it is given by the formula 
$$
{\mathcal I}_{s}(M) = 32\pi^{2}c_{1}^{2}(X),
$$
where $X$ is the minimal model of $M$, in the sense of 
Kodaira. 

One main purpose of the present article is to
similarly calculate ${\mathcal I}_{r}$ for the complex surfaces. 
Unlike the ${\cal I}_{s}$, this invariant is changed by
blowing up, and in Theorem \ref{l2ric} we will see that
it is given by 
$$
{\mathcal I}_{r}(M) = 8\pi^{2}[c_{1}^{2}(X)+ k]
$$
where $k$ is the number of points of $X$ which must 
be blown up in order to obtain $M$. 

This, however, is just one application of some new 
curvature estimates we will develop here. Fundamentally,
the  story  is not  
about  
 the Ricci curvature at all, but instead principally concerns  the   
r\^ole of  the {\em  self-dual Weyl curvature}.
Recall that the $2$-forms on an oriented $4$-manifold
decompose as 
$$\Lambda^{2}= \Lambda^{+}\oplus \Lambda^{-},$$ 
where $\Lambda^{\pm}$ is the $(\pm 1)$ eigenspace
of  Hodge  star operator $\star$. Thinking of the 
curvature tensor $\mathcal R$ as a linear map
$\Lambda^{2}\to \Lambda^{2}$, we thus get 
a decomposition \cite{st}
$$
{\mathcal R}=
\left(
\mbox{
\begin{tabular}{c|c}
&\\
$W_++\frac{s}{12}$&$B$\\ &\\
\cline{1-2}&\\
$B^{*}$ & $W_{-}+\frac{s}{12}$\\&\\
\end{tabular}
} \right) 
$$
into irreducible pieces. 
Here the
self-dual and anti-self-dual Weyl curvatures 
$W_\pm$ are the trace-free pieces of the appropriate blocks.
The scalar curvature  $s$ is understood to act by scalar multiplication,
and $B$ amounts to the trace-free part 
$\stackrel{\circ}{r}=r-\frac{s}{4}g$  of the Ricci curvature, 
acting on anti-self-dual 2-forms by
$$\psi_{ab} \mapsto ~ 
 \stackrel{\circ}{r}_{ac}{\psi^c}_{b}-
\stackrel{\circ}{r}_{bc}{\psi^c}_{a}.$$

The key estimate of the article is to be found in 
Theorem \ref{L2}: if $(M,g)$ is an oriented Riemannian 
$4$-manifold, and if, for a fixed spin$^{c}$ structure on $M$, 
the Seiberg-Witten equations have a solution for every metric conformally
related to $g$,  then the curvature of $g$ satisfies
$$	\int_{M}\left({\frac{2}{3}s-2\sqrt{\frac{2}{3}}|W_{+}|}\right)^{2} 
	d\mu 
\geq 32\pi^{2}(c_{1}^{+})^{2},$$
where $c_{1}^{+}$ is the self-dual part, with respect to $g$,  of 
the first Chern class of the spin$^{c}$ structure.
Information about the Ricci curvature is then indirectly
extracted from this, by means of the Gauss-Bonnet-like 
formula  
\begin{equation}
(2\chi + 3\tau ) (M) 
=\frac{1}{4\pi^{2}}\int_{M}\left(
2|W_{+}|^{2}+\frac{s^{2}}{24} -\frac{|\stackrel{\circ}{r}|^{2}}{2}
\right) d\mu ,	
	\label{bonnet}
\end{equation}
which holds for 
 every Riemannian metric $g$  on $M$.
Here $\chi(M)$ and $\tau(M)$ denote the Euler characteristic and signature
of $M$,
respectively. 

The above estimate
 should be compared and contrasted to the 
analogous inequality 
$$
\int_{M}s^{2} 
	d\mu 
\geq 32\pi^{2}(c_{1}^{+})^{2}.$$
for the scalar curvature, keeping in mind  that 
$|s|=2\sqrt{6}|W_{+}|$ for any K\"ahler metric.
This seems all the more remarkable insofar as 
the older estimate is only  saturated by 
constant-scalar-curvature 
K\"ahler metrics, while  the new estimate is  
apparently saturated by a  larger class  of
{\em almost-K\"ahler} metrics. This reflects
an under-utilized aspect of  Taubes' construction of solutions of the 
Seiberg-Witten equations on  symplectic
manifolds \cite{taubes,t2}, and  would appear 
to be a promising avenue for   further
research.  

 These estimates also give one new obstructions to  the existence 
 for Einstein metrics. Recall that that a smooth Riemannian metric 
 $g$ is said to be {\em  Einstein} 
if its Ricci curvature $r$ is a constant multiple of 
the metric: 
$$r=\lambda g .$$
Not every smooth compact oriented 4-manifold $M$
admits such a metric. Indeed, 
a well-known necessary condition  is that 
$M$ must satisfy \cite{hit,thorpe,bes} 
the 
Hitchin-Thorpe inequality $2\chi  (M) \geq  3 |\tau (M)|$,
where again $\chi$ and $\tau$ denote the signature and
Euler characteristic. Indeed, this is an immediate consequence of 
(\ref{bonnet}), 
since the Einstein condition may be rewritten as 
$\stackrel{\circ}{r}= 0$, and $\stackrel{\circ}{r}$
 makes the only negative contribution to the integrand. 
Notice, however, that one could strengthen the conclusion 
given an {\em a priori} lower bound 
for the  term
$$\frac{1}{4\pi^{2}}\int_{M}\left(
2|W_{+}|^{2}+\frac{s^{2}}{24} 
\right) d\mu .$$
But a lower bound is not difficult to 
extract from the above estimate, and this allows one     
 to show
(Corollary \ref{blowup}) that if a minimal
complex surface $X$ is blown up $k$ times,
the resulting complex surface
$$M= X \# k \overline{\bcp}_{2}$$
does not admit Einstein metrics if $k\geq \frac{1}{3}c_{1}^{2}(X)$.
This improves previous results \cite{lno,lebweyl},
where the same conclusion is reached for larger values of 
$k$.

The invariant ${\mathcal I}_{\mathcal R}$ introduced at the 
beginning of this introduction is often easy to compute in 
dimension $4$, because \cite{bes}
$${\mathcal I}_{\mathcal R}(M)\geq 8\pi^{2}\chi (M),$$
with equality if $M$ admits an Einstein metric. In particular, 
it is trivial to read off this invariant 
for a compact hyperbolic $4$-manifold ${\mathcal H}^{4}/\Gamma$
or for a  complex-hyperbolic $4$-manifold ${\mathbb C}
{\mathcal H}_{2}/\Gamma$. (The complex hyperbolic plane
${\mathbb C}
{\mathcal H}_{2}$ may abstractly be defined as the
symmetric space $SU(2,1)/U(2)$, but it is typically  more
useful to think of it as 
the unit ball in ${\mathbb C}^{2}$, equipped with the
Bergmann metric.) 
It might be tempting to assume that this tells one
everything there is to know about the sectional curvatures
of metrics on these spaces; but in reality, that  is wide  
 off the mark!  
For example, one might ask whether the standard
metric, suitably normalized, has least volume
among all metrics of sectional curvature
$K\geq -1$. For  real-hyperbolic manifolds
the answer is affirmative, and indeed this  holds in all
dimensions; however, the proof, due to Besson-Courtois-Gallot 
\cite{bcg2},
depends not only on a remarkable inequality concerning volume entropy,
but also on an optimal application of Bishop's inequality. 
For complex-hyperbolic manifolds, the latter breaks down, 
and the question  therefore remains open.

While the above question cannot be settled by present means,  our
estimates do allow one to answer a closely  
related question. 
Notice that 
a $4$-manifold with $K\geq -1$ also satisfies 
$$\frac{1}{2}\left(K+\frac{s}{12}\right)\geq -1.$$ 
Given a complex hyperbolic $4$-manifold, we are
able to show (Theorem \ref{minvol}) that, among metrics on 
a complex-hyperbolic manifold satisfying this curvature constraint, 
 an appropriate  multiple  
of the 
 standard metric has least volume. Similar
results also hold  metrics 
subject to the curvature constraint
$$
tK+(1-t)\frac{s}{12}\geq -1 
$$
for any  constant $t\in [0,1/2]$.

\section{Weyl  Curvature Estimates}
\label{weylest} 
In this section, we derive new  
$L^{2}$ estimates for combinations of
 the  Weyl and scalar  curvatures of certain  Riemannian 
$4$-manifolds. These considerably refine 
the  estimates previously found in \cite{lebweyl}.

Let $M$ be a smooth, compact, oriented 4-manifold. Each 
Riemannian metric $g$ on $M$ then determines
a direct sum decomposition
$$H^{2}(M, {\Bbb R}) =
{\cal H}^{+}_{g}\oplus {\cal H}^{-}_{g},$$
where ${\cal H}^{+}_{g}$  (respectively, ${\cal H}^{-}_{g}$)
 consists of those cohomology
classes for which the   harmonic  representative  is 
 self-dual  (respectively,  anti-self-dual).
 The non-negative integer
$b_{+}(M)= \dim {\cal H}^{+}_{g}$  is independent of 
$g$, and we will henceforth always assume it to be positive. 
It is thus  natural to consider the set of metrics $g$
for which ${\cal H}^{+}_{g}=H$ for some 
fixed $b_{+}(M)$-dimensional subspace $H\subset H^{2}(M, {\Bbb R})$;
such metrics will be said to be {\em $H$-adapted}. 
 Assuming there is at least one $H$-adapted metric, we will then 
 say that $H$ is a {\em polarization} of $M$, and 
 call the pair 
$(M,H)$  a {\em polarized 4-manifold} \cite{lpm}. Notice that
the restriction of
the intersection pairing
$$\smile : H^{2}(M, {\Bbb R})\times H^{2}(M, {\Bbb R})\to {\Bbb R}$$
to $H$ 
is then positive definite, and that $H\subset H^{2}$
is   maximal among  subspaces with this property.

Let $c$ be a spin$^{c}$ structure on $M$. Then $c$ determines 
a Hermitian line-bundle $L\to M$ with 
$$c_{1}(L)\equiv w_{2}(M)\bmod 2,$$
and for each metric $g$ we also  have  
rank-2 complex vector bundles
${\mathbb V}_{\pm}\to M$ which formally satisfy 
$${\mathbb V}_{\pm}={\Bbb S}_{\pm}\otimes L^{1/2},$$
where ${\Bbb S}_{\pm}$ are the locally-defined left- and right-handed 
spinor bundles of $g$. Given a polarization $H$ on $M$, we 
will then use $c_{1}^{+}$ to denote the 
orthogonal projection of $c_{1}(L)$ into $H$ with 
respect to the intersection form. If $g$ is a particular
metric with ${\cal H}^{+}_{g}=H$, we will also freely use
$c_{1}(L)$ to denote the $g$-harmonic 2-form representing the 
corresponding 
de Rham class, and use $c_{1}^{+}$ to denote  its self-dual part. 
For example, if a choice of $H$-compatible metric $g$
has already been made,  the number 
$$|c_{1}^{+}|:=\sqrt{(c_{1}^{+})^{2}}$$
may freely be identified with the
L$^{2}$-norm of the self-dual $g$-harmonic form 
denoted by $c_{1}^{+}$.

 For each Riemannian metric $g$,  
the Seiberg-Witten equations \cite{witten} 
\begin{eqnarray} D_{A}\Phi &=&0\label{drc}\\
 F_{A}^+&=&i \sigma(\Phi)\label{sd}\end{eqnarray}
are 
equations 
for an unknown Hermitian connection $A$ on $L$
and an unknown 
smooth section $\Phi$ of ${\mathbb V}_+$.
Here $D_{A}$ is the Dirac operator coupled to $A$, 
and $\sigma : 
{\mathbb V}_{+}\to \Lambda^{+}$ is a certain 
canonical real-quadratic map. The latter formally  arises from the fact that
${\mathbb V}_{+}={\Bbb S}_{+}\otimes L^{1/2}$, whereas 
$\Lambda^{+}\otimes {\mathbb C}= \odot^{2}{\Bbb S}_{+}$;
it 
is invariant under parallel transport, and 
satisfies $|\sigma (\Phi )|^{2}= |\Phi |^{4}/8$.
In conjunction with (\ref{sd}), 
the latter immediately gives us the important inequality 
\begin{equation}
	\int_{M}|\Phi |^{4}d\mu \geq 32\pi^{2}(c_{1}^{+})^{2}
	\label{harm}
\end{equation}
because $2\pi c_{1}^{+}$ is the harmonic part of 
$-\sigma (\Phi )$.

In this paper, we will primarily  be interested in $4$-manifolds for which 
there is a solution of the Seiberg-Witten equations 
for each metric. Let us make this more precise by introducing some
 terminology; cf. \cite{K}.

\begin{defn}
Let $M$ be a smooth compact oriented $4$-manifold,
let $H\subset H^{2}(M, {\mathbb R})$ be a polarization of $M$, 
and let $c$ be a spin$^{c}$ structure on $M$. Then we will say that 
{\em $c$ is a monopole class for $(M,H)$} if the Seiberg-Witten 
equations (\ref{drc}--\ref{sd})
have a solution for every $H$-adapted metric $g$. 
\end{defn}

This definition will be useful in practice, of course, only
because of the existence of Seiberg-Witten invariants \cite{KM,witten} . 
For example, if  a spin$^c$ structure satisfies 
 $[c_{1}(L)]^{2}=(2\chi + 3\tau )(M)$,  and if 
  $c_{1}^{+}\neq 0$ relative to the polarization $H={\cal H}^{+}_{g}$, 
 then   the  {\em Seiberg-Witten invariant}   
${\mathcal S \mathcal W}_c(M ; H)$ can be defined as the number of solutions,
modulo gauge transformations and  counted with orientations, of a  
generic perturbation 
\begin{eqnarray*} D_{A}\Phi &=&0  \\
 iF^+_A+\sigma (\Phi ) &=& \phi  \end{eqnarray*}
 of  (\ref{drc}--\ref{sd}), 
where $\phi$ is a smooth self-dual 2-form of small 
$L^{2}$ norm. If this invariant is non-zero, then $c$ is a 
monopole class for $(M,H)$, where $H$ is the polarization
determined by the metric $g$. If $b^{+}(M)\geq 2$, 
 ${\mathcal SW}_{c}(M,H)$ is actually independent of
 the polarization $H$; when $b^{+}(M)=1$, by contrast,
 it is well defined only for those polarizations for
 which $c_{1}^{+}\neq 0$, and its value typically depends on whether
 $c_{1}^{+}$ is a future-pointing or past-pointing 
 time-like vector in the Lorentzian vector space $H^{2}(M, {\mathbb R})$.

More generally \cite{ozsz,lebyond},
suppose that we have a spin$^{c}$ structure such that
$$\ell = \frac{[c_{1}(L)]^{2}-(2\chi + 3\tau )(M)}{4}$$
is a non-negative integer. Then the moduli space ${\mathcal M}_{c,g}$
of  
gauge-equivalence classes of solutions 
of a generic perturbation of  (\ref{drc}--\ref{sd}) 
 is a
smooth compact $\ell$-manifold; moreover, it acquires a 
canonical orientation once we choose an orientation 
 for the vector space $H^{1}(M,{\mathbb 
R})\oplus H$. 
Fix $\ell$ loops $\beta_{1}, \ldots  , 
 \beta_{\ell}$ 
 in $M$, and define a smooth
 map 
 ${\mathcal M}_{c, g }\to T^{\ell}$ from the 
 moduli space to the $\ell$-torus 
 by sending the gauge-equivalence class $[(\Phi, A)]$ to 
 the holonomies of the $U(1)$ connection $A$ around 
 the $\ell$ given loops. 
 The homotopy class  of this map only depends on
 the homology classes $[\beta_{i}]\in H_{1}(M, {\mathbb Z})$,
 one we may therefore define 
 ${\mathcal S \mathcal W}_{c}(M, [\beta_{1}], \cdots , 
 [\beta_{\ell}] ; H)\in {\mathbb Z}$ to be the degree  of this 
 map. Again, if this invariant is non-zero, $c$ is a monopole class
 of $(M,H)$.

 Many of the most remarkable consequences of  Seiberg-Witten 
 theory stem from the fact that the equations (\ref{drc}--\ref{sd})  
imply the Weitzenb\"ock formula
\begin{equation}\label{wb}
0= 4\nabla^*\nabla \Phi +  s \Phi +|\Phi|^2\Phi ,
\end{equation}
where $s$ denotes the scalar curvature of $g$.
Taking the inner product with $\Phi$, it follows that
\begin{equation}
	0=2\Delta |\Phi |^{2} +4|\nabla \Phi |^{2} + s |\Phi |^{2}+ |\Phi 
|^{4} .
	\label{zero}
\end{equation}
If we  multiply (\ref{zero}) by 
 $|\Phi |^{2}$ and integrate,  we  have
$$0 = \int_{M}\left[ 2\left| ~d|\Phi |^{2} ~\right|^{2} + 4 |\Phi |^{2}|\nabla 
\Phi |^{2}+ s|\Phi |^{4}+ |\Phi |^{6}
\right] d\mu_{g} ,$$
so that 
\begin{equation}
	\int (-s) |\Phi |^{4} d\mu \geq 4 \int |\Phi |^{2} |\nabla \Phi 
|^{2} d\mu  + \int  |\Phi |^{6} d\mu .
	\label{two}
\end{equation}
We will  see in a moment  
that this implies some remarkable estimates
 for the Weyl curvature of suitable $4$-manifolds. 

Before doing so, however, let us first recall that the self-dual 
  Weyl tensor at a point $x$ of an
   oriented Riemannian $4$-manifold $(M,g)$ may be
  viewed as a trace-free endomorphism $W_{+}(x): \Lambda^{+}_{x}
  \to \Lambda^{+}_{x}$ of the  self-dual $2$-forms at $x$. 
 We will use $w(x)$ denotes its lowest eigenvalue; notice that this is 
 automatically a Lipschitz continuous function $w: M\to (-\infty, 0]$. 
  Let us also fix the notational convention that, for any
  real-valued function $f: M\to {\mathbb R}$, 
  $f_{-}: M\to (-\infty, 0]$ is defined by
  $f_{-}(x)= \min (f(x), 0)$. With these preliminaries, 
  we are now prepared to prove our first, crucial result:

 \begin{prop}\label{L3} 
 Let $(M,g)$ be a compact oriented Riemannian $4$-manifold
  on which there is a solution of the Seiberg-Witten 
  equations. Let 
 $c_{1}(L)$ be the  first Chern class of the relevant spin$^{c}$ structure,
 and let $c_{1}^{+}$ denote its self-dual part with respect to $g$.
 Then 
 \begin{equation}
 	V^{1/3}\left( \int_{M}|(\frac{2}{3}s_{g}+ {2}w_{g} )_{-}|^{3}d\mu 
 \right)^{2/3}\geq 32\pi^{2} (c_{1}^{+})^{2},
 	\label{crux}
 \end{equation}
 where $V=\mv (M,g)= \int_{M}d\mu_{g}$ is the total volume of $(M,g)$. 
 \end{prop}
 \begin{proof} 
 Any self-dual  2-form $\psi$ on any oriented 4-manifold satisfies 
 the 
Weitzenb\"ock formula \cite{bourg}
$$(d+d^{*})^{2}\psi = \nabla^{*}\nabla \psi - 2W_{+}(\psi , 
\cdot ) + \frac{s}{3} \psi,$$
where $W_{+}$ is the self-dual Weyl tensor. 
It follows that 
$$
\int_{M}(-2W_{+})(\psi , \psi ) \geq 
\int_{M}(-\frac{s}{3})|\psi |^{2}~d\mu -
 \int_{M} |\nabla \psi |^{2}
 ~d\mu , 	
$$
so that 
 $$
  -\int_{M}2w|\psi |^{2} \geq 
\int_{M}(-\frac{s}{3})|\psi |^{2}~d\mu -
\int_{M} |\nabla \psi |^{2}
 ~d\mu , 
  $$ 
  and hence
   $$
  -\int_{M}(\frac{2}{3}s+2w)|\psi |^{2} \geq 
\int_{M}(-s)|\psi |^{2}~d\mu -
\int_{M} |\nabla \psi |^{2}
 ~d\mu . 
  $$ 
On the other hand, the  particular self-dual 2-form $\varphi = \sigma (\Phi 
)=-iF_{A}^{+}$ satisfies 
\begin{eqnarray*}
	|\varphi |^{2} & = & \frac{1}{8}|\Phi |^{4}  ,\\
	|\nabla \varphi |^{2} & \leq  & \frac{1}{2} |\Phi |^{2}|\nabla \Phi |^{2} .
\end{eqnarray*}
 Setting $\psi = \varphi$,  
  we thus have
  $$
   -\int_{M}(\frac{2}{3}s+2w)|\Phi |^{4} \geq 
\int_{M}(- s)|\Phi |^{4}~d\mu -
 4\int_{M} |\Phi |^{2} |\nabla \Phi |^{2}
 ~d\mu .
  $$
 But   (\ref{two})
 tells us that 
 $$\int_{M}(- s)|\Phi |^{4}~d\mu -
 4\int_{M} |\Phi |^{2} |\nabla \Phi |^{2}
 ~d\mu \geq \int_{M}|\Phi |^{6}~d\mu ,$$
 so we obtain 
 $$-\int_{M}(\frac{2}{3}s+2w)_{-}|\Phi |^{4} d\mu \geq
 -\int_{M}(\frac{2}{3}s+2w)|\Phi |^{4} d\mu \geq  
 \int_{M}|\Phi |^{6}~d\mu .$$
By the H\"older inequality, we thus have
$$\left( \int |(\frac{2}{3}s+2w)_{-}|^{3}d\mu \right)^{1/3}
\left( \int |\Phi |^{6} d\mu \right)^{2/3}
\geq  
 \int|\Phi |^{6}~d\mu ,$$
 Since  the H\"older inequality also tells us 
 that
 $$
 \int|\Phi |^{6}~d\mu  \geq V^{-1/2}\left(\int |\Phi |^{4}d\mu \right)^{3/2} ,
 $$
 we thus have 
$$ V^{1/3}\left( \int_{M}|(\frac{2}{3}s_{g}+ {2}w_{g} )_{-}|^{3}d\mu 
 \right)^{2/3}\geq 
 \int |\Phi |^{4}d\mu
 \geq 32\pi^{2} (c_{1}^{+})^{2},$$
 as claimed. 
 \end{proof}
 
 While this result is the key to 
  everything that follows, its direct utility is  limited by 
 the fact that it is an $L^{3}$, rather than an $L^{2}$, estimate. 
 Fortunately, however, we will be able to extract an  $L^{2}$ estimate
 by means of a conformal rescaling  trick, the general idea of which 
 is drawn from    Gursky \cite{gu}:

 \begin{lem}\label{gursky}
 Let $(M, \gamma)$ be a compact oriented $4$-manifold
 with a fixed smooth conformal class of Riemannian metrics. 
 Suppose, moreover,  that $\gamma$ does not contain a
 metric
 of positive scalar curvature. Then, for any $\alpha \in (0,1)$,
  there is a  metric $g_{\gamma}\in \gamma$ of differentiability class
 $C^{2,\alpha}$ for which 
 $s+ 3w$ is a non-positive constant. 
 \end{lem}
 \begin{proof} Let $g_{0}\in \gamma$ be a fixed smooth back-ground metric,
 and notice that, for $g= u^{2}g_{0}\in \gamma$, 
  the 
 function ${\mathfrak S}_{g} = s_{g}+ 3w_{g}$ is 
 given by 
 $${\mathfrak S}_{g} = u^{-3}\left( 6 \Delta_{g_{0}}u + {\mathfrak 
 S}_{g_{0}}u \right)$$
 by virtue of 
 the weighted conformal invariance of $W_{+}$. Thus 
  $\mathfrak S$ behaves under conformal rescaling
 just like the scalar curvature, despite the fact that it might well
 only be a Lipschitz function. In order to find a suitable choice of 
$u$, we  therefore attempt to minimize 
$${\mathcal F}(u) = \frac{\int_{M}\left(6|du|^{2}+ {\mathfrak 
 S}_{g_{0}}u^{2}\right)d\mu_{g_{0}}}{\left(
 \int_{M}u^{4}d\mu_{g_{0}}\right)^{1/2}}
$$
on the positive sector of the unit sphere 
in the Sobolev space $L^{2}_{1}(M,g_{0})$. 
Yamabe's ansatz  for 
doing this \cite{yam,aubin0} is to minimize the functionals
$${\mathcal F}_{\epsilon}(u) = \frac{\int_{M}\left(6|du|^{2}+ {\mathfrak 
 S}_{g_{0}}u^{2}\right)d\mu_{g_{0}}}{\left(
 \int_{M}u^{4-\epsilon}d\mu_{g_{0}}\right)^{1/(2-\frac{\epsilon}{2})}}
$$
and then take the limit of the minimizers as $\epsilon \searrow 0$.
For each $\epsilon > 0$, the existence of  minimizers follows from 
 the Sobolev
embedding theorem; indeed, following 
the proof of 
 \cite[Theomem 5.5]{aubin0}, 
but deleting the very last sentence,  
one  obtains a   positive function $u_{\epsilon}$ of class 
$C^{2, \alpha}$
which solves
 \begin{equation}
 6 \Delta_{g_{0}}u + {\mathfrak 
 S}_{g_{0}}u= c_{\epsilon} u^{3-\epsilon} ,	
 	\label{el}
 \end{equation} 
 where $c_{\epsilon}$ is the infimum of ${\mathcal F}_{\epsilon}$
 on the positive sector of the 
 unit sphere of $L^{2}_{1}(M,g_{0})$. 
   The convergence as $\epsilon \searrow 0$ 
 then follows from an observation of  
 Trudinger \cite{trud}. Namely, since 
 ${\mathfrak S} \leq s$ for any metric, 
 and since there is nothing to prove if ${\mathfrak S} = s$,
  our hypothesis that there is no metric
 of positive scalar curvature in $\gamma$ allows us to assume that 
 that $c_{0}=\inf c_{\epsilon}
 < 0$. Inspection of (\ref{el}) at a maximum then gives us 
 the $C^{0}$ estimate
 $$u_{\epsilon} \leq 1+\sqrt{\left|\frac{\inf {\mathfrak 
 S}_{g_{0}}}{c_{0}}\right|} 
 $$
 for all small $\varepsilon$.
 This implies \cite[Theorem 6.5]{aubin0} that   
 $u=\lim_{\epsilon\to 0} u_{\epsilon}$ exists in 
 $C^{1}$, and is a weak solution of 
 $$ 6 \Delta_{g_{0}}u + {\mathfrak 
 S}_{g_{0}}u= c_{0} u^{3} .$$
 Schauder theory then tells us that $u$ is actually a
 $C^{2,\alpha}$ function. In particular, $g_{\gamma}=u^{2}g_{0}$ 
 is a $C^{2,\alpha}$ metric, with curvature  satisfying 
  ${\mathfrak S}_{g_{\gamma}}=c_{0}$
  in the classical sense.  
 \end{proof}

\begin{thm} \label{L2}
Let $(M,H)$ be a  polarized smooth compact  
oriented  4-manifold, and let $c$ be a monopole class for 
$(M,H)$. Let $c_1(L)\in H^2(M, {\Bbb R})$ denote 
the anti-canonical class of $c$, and 
  let $c_1^+\neq 0$ be its orthogonal projection to
$H$ with respect to the intersection form. 
Then  every  
$H$-adapted Riemannian metric $g$ 
satisfies 
\begin{equation}
	\int_{M}\left({\frac{2}{3}s-2\sqrt{\frac{2}{3}}|W_{+}|}\right)^{2} 
	d\mu 
\geq 32\pi^{2}(c_{1}^{+})^{2}.
	\label{central}
\end{equation}
\end{thm}
\begin{proof} Let $\gamma$ be the conformal class of some
$H$-adapted metric $g$. Since the Hodge star operator is 
conformally invariant in the middle dimension, every metric
in $\gamma$ is also $H$-adapted. Since $c$ is a monopole class for
$H$, it
therefore follows that $\gamma$  does
not contain any metrics of positive scalar curvature. 
Lemma \ref{gursky} therefore applies, and tells us that 
the conformal class $\gamma$ contains 
a metric $g_{\gamma}$ for which that $\frac{2}{3}s+2w= \frac{2}{3}{\mathfrak S}$
is a non-positive constant. For this metric, 
one then has
$$\int_{M}|(\frac{2}{3}s_{g_{\gamma}}+2w_{g_{\gamma}}
 )_{-}|^{2}d\mu_{g_{\gamma}} =
V^{1/3}_{g_{\gamma}}\left( \int_{M}|(\frac{2}{3}s_{g_{\gamma}}+ 
{2}w_{g_{\gamma}} )_{-}|^{3}d\mu_{g_{\gamma}} 
 \right)^{2/3},$$
 so that
 $$
 \int_{M}|(\frac{2}{3}s_{g_{\gamma}}+2w_{g_{\gamma}}
 )_{-}|^{2}d\mu_{g_{\gamma}} \geq 
  32\pi^{2} (c_{1}^{+})^{2} 
 $$
 by Proposition \ref{L3}. 
 Thus we at least have an  $L^{2}$ estimate concerning the 
 conformally related metric $g_{\gamma}$.
 
Let us now compare the left-hand side with analogous 
expression for 
 the given metric $g$. To do so, 
we express $g$ in the form $g=u^{2}g_{\gamma}$,
where $u$ is a positive $C^{2}$ function, and observe that 
\begin{eqnarray*}
	\int_{M} (\frac{2}{3}s_{g}+2w_{g}
 )_{-} u^{2}d\mu_{g_{\gamma}} &\leq &
	\int_{M} \left(\frac{2}{3}s_{g}+2w_{g}
 \right) u^{2}d\mu_{g_{\gamma}}\\ & = & \frac{2}{3}\int 
 {\mathfrak S}_{g}u^{2}d\mu_{g_{\gamma}}  \\
	 & = & \frac{2}{3}\int u^{-3}\left(6\Delta_{g_{\gamma}}u+
	 {\mathfrak S}_{g_{\gamma}}u\right) u^{2}d\mu_{g_{\gamma}}
	 \\&=& \frac{2}{3}\int \left(-6 u^{-2} |du|^{2}_{g_{\gamma}}+
	 {\mathfrak S}_{g_{\gamma}}\right)  d\mu_{g_{\gamma}}
	 \\ &\leq & \frac{2}{3}\int 
	 {\mathfrak S}_{g_{\gamma}}  d\mu_{g_{\gamma}}\\
	 &=& \int_{M} (\frac{2}{3}s_{g_{\gamma}}+2w_{g_{\gamma}})_{-} 
	 d\mu_{g_{\gamma}} .
\end{eqnarray*}
Applying Cauchy-Schwarz, we thus have
\begin{eqnarray*}
 -V^{1/2}_{g_{\gamma}}\left(\int |(\frac{2}{3}s_{g}+2w_{g}
 )_{-}|^{2} d\mu_{g}\right)^{1/2}	 & = &
	 -V^{1/2}_{g_{\gamma}}\left(\int \left(\frac{2}{3}s_{g}+2w_{g}
 \right)^{2} u^{4}d\mu_{g_{\gamma}}\right)^{1/2}\\
 & \leq  & 
 \int_{M} \left(\frac{2}{3}s_{g}+2w_{g}
 \right) u^{2}d\mu_{g_{\gamma}}
 \\&\leq &\int_{M} \left(\frac{2}{3}s_{g_{\gamma}}
 +2w_{g_{\gamma}}\right) 
	 d\mu_{g_{\gamma}}  \\
	 & = &  -
	 V^{1/2}_{g_{\gamma}}\left(\int 
	 |(\frac{2}{3}s_{g_{\gamma}}+2w_{g_{\gamma}}
	  )_{-}|^{2}d\mu_{g_{\gamma}}\right)^{1/2} , 
\end{eqnarray*}
and hence  
$$\int_{M}|(\frac{2}{3}s_{g}+2w_{g} )_{-}|^{2}d\mu_{g} \geq 
\int_{M}|(\frac{2}{3}s_{g_{\gamma}}+2w_{g_{\gamma}}
 )_{-}|^{2}d\mu_{g_{\gamma}} .
$$
This shows that 
\begin{equation}
	\int_{M}|(\frac{2}{3}s_{g}+2w_{g})_{-} |^{2}d\mu_{g} 
	\geq 32\pi^{2}(c_{1}^{+})^{2}
	\label{yes}
\end{equation}
for every $H$-adapted metric $g$. 

Finally, we observe that 
$$ -\sqrt{\frac{2}{3}}|W_{+}|_{g}\leq w_{g}$$
for any metric, simply because $W_{+}$ is trace-free.
It follows that 
$$\left(\frac{2}{3}s-2\sqrt{\frac{2}{3}}|W_{+}|_{g} \right)^{2}
\geq |(\frac{2}{3}s_{g}+2w_{g})_{-} |^{2}$$
at every point, and (\ref{yes}) therefore implies the
desired inequality (\ref{central}). 
  \end{proof}

The inequalities we have just derived are certainly sharp, in the 
sense that equality is attained by any K\"ahler metric of constant
negative scalar curvature. However, we will now see that 
equality also holds in principle for a broader class of metrics. To make this 
precise, suppose that $(M,\omega )$ is a symplectic $4$-manifold, 
and let $g$ be a Riemannian metric such that 
$$\omega (\cdot , \cdot ) = g (J \cdot , \cdot )$$
for some almost complex structure $J$ on $M$; 
equivalently, suppose that $\omega$ is self-dual with respect to 
$g$ and has constant norm $|\omega |_{g}= \sqrt{2}$. 
One then says that $g$ is an {\em almost-K\"ahler metric},
with {\em almost-K\"ahler form} $\omega$. Evidently, 
$g$ is actually K\"ahler if $J$ is integrable, but Gromov's
theory of pseudo-holomorphic curves \cite{gromsym} and 
Taubes' characterization
\cite{taubes,t2} of the 
Seiberg-Witten invariants of symplectic $4$-manifolds 
have  conclusively swept the study of such metrics into
the mathematical main-stream.

In order to state our next result, we will need  
 extra  terminology in this connection.  The $\ast$-scalar curvature of 
a $4$-dimensional almost-K\"ahler manifold $(M,g,\omega )$ 
is by definition the function 
$$s^{*}= \frac{s}{3}+2W_{+}(\omega , \omega ) .$$
This definition is motivated  by the fact that 
the scalar curvature $s$ and  $\ast$-scalar curvature
$s^{*}$ coincide for any K\"ahler manifold. Notice that
we automatically have the inequality 
$$s^{*}\geq \frac{s}{3} + 4w ,$$
with equality iff the almost-K\"ahler form $\omega$
belongs to the lowest eigenspace of $W_{+}$.

\begin{prop}\label{aka} Suppose that equality holds in either (\ref{crux})
or (\ref{yes}).
Then   $g$ is an almost-K\"ahler metric
 with the following properties:
 \begin{itemize}
 \item the almost-K\"ahler form $\omega$ belongs, at each $x\in M$,
  to the 
 lowest eigenspace of $W_{+}$; 
 \item the sum $s+s^{*}$ of  the  scalar and  $\ast$-scalar curvatures
 is a non-positive {constant};
 \item the  first Chern class  of $(M,\omega )$ is of
 type $(1,1)$, and coincides with  first Chern class $c_{1}(L)$
 of the relevant spin$^{c}$ structure. 
 \end{itemize}
 If equality holds in (\ref{central}),  these same conclusions hold,
 but, in addition,  
 \begin{itemize}
 \item the self-dual Weyl curvature $W_{+}$ 
 is invariant with respect to the 
 almost-complex structure $J$ of $(M,g,\omega )$. 
 \end{itemize} 
 Conversely,
  any compact almost-K\"ahler $4$-manifold satisfying these
 conditions saturates the corresponding inequalities, and 
 its first Chern class is a monopole class of the corresponding 
 polarization. 
\end{prop}
\begin{proof}
If equality	holds in (\ref{crux}), the last inequality in  the proof 
of Proposition \ref{L3}  forces  
 $\sigma (\Phi )$ to be	
  	the	harmonic representative	of $-2\pi c_{1}^{+}$. 
  Moreover,   the H\"older inequalities used in that proof shows
  	that  equality can only	hold if	
  	$|\Phi|^{2}$ and $-(\frac{2}{3}s+2w)$ are  constant, 
  	and moreover	equal.
  	Thus the harmonic form $\varphi= \sigma (\Phi )$	has	constant norm, 
  	so that $\omega = -\sqrt{2} \varphi /|\varphi |$ is an
  	almost-K\"ahler form compatible with $g$. Moreover, 
  	$L=\det ({\mathbb V}_{+})= {\mathbb V}_{+}/\mbox{span} (\Phi )$  coincides 
  	with the anti-canonical 
  	line bundle associated with the almost-complex structure $J$. 
  	Since 
  	$c_{1}^{+}$ is now a (non-positive) multiple of $[\omega ]$
  	and since the volume $V=\mv (M,g)$ of an almost-K\"ahler $4$-manifold is 
  	exactly $[\omega ]^{2}/2$, the fact that 
  	 $\frac{2}{3}s+2w$ is a non-positive constant tells us that 
  	 we now have 
  	$$\int (\frac{2}{3}s+2w) d\mu = -V^{1/2}4\sqrt{2}\pi |c_{1}^{+}| 
  	= 
  	4\pi c_{1}(L)\cdot [\omega ],$$
  	since by assumption equality holds in (\ref{crux}). It thus follows 
  	that
  	$$\int \frac{1}{2}(s+s^{*}) d\mu \geq 4\pi c_{1}(L)\cdot 
  	[\omega ],$$
  	with equality iff $\omega$ everywhere belongs to the lowest
  	eigenspace of $W_{+}$. 
  	However, Blair \cite{blair} has shown that 
  	\begin{equation}
  		\int \frac{1}{2}(s+s^{*}) d\mu= 4\pi c_{1}\cdot 
  	[\omega ]
  		\label{blr}
  	\end{equation}
  	for any almost-K\"ahler $4$-manifold. We thus
  	conclude that $\omega$  is everywhere in the lowest eigenspace of 
  	$W_{+}$. 
  	
  	If we instead have equality in (\ref{yes}) or (\ref{central}),
  	we must in particular be at the minimum of 
  	$$\int |(\frac{2}{3}s + 2w)_{-}|^{2}d\mu$$
  	among metrics in the given conformal class. This then implies that 
  	$\frac{2}{3}s + 2w$ is constant, and we therefore also have
  	equality in (\ref{crux}), and the claim follows from the 
  	previous argument. 
  	
  	The converse assertions  follow from 
  	\cite{taubes} and 
  	 (\ref{blr}). Details
  	are left to the interested reader.
  \end{proof}
  
  In a na\"{\i}ve sense, this provides a complete characterization of
  the metrics which saturate (\ref{crux}), (\ref{central}), and 
  (\ref{yes}). It remains to be seen, however, 
  whether such  metrics can ever be {\em strictly} almost-K\"ahler,
  in the sense that the almost-complex structure $J$ fails to be
  integrable. For example, a beautiful recent result of J. Armstrong 
   \cite{rmstg} asserts the non-existence of compact, Einstein,  
  strictly almost-K\"ahler 
  $4$-manifolds on which  $\omega$ is everywhere
  an eigenvector of $W_{+}$. Thus:

  \begin{cor}
  Suppose that $g$ is an Einstein metric 
   which saturates
 (\ref{crux}), (\ref{central}), or  
  (\ref{yes}). Then $g$ is actually K\"ahler-Einstein. 
  \end{cor}
  
  Similarly, a recent result of Apostolov-Armstrong-Dr{\u{a}}ghici \cite{ad} 
  implies that an 
  almost-K\"ahler metric
  saturating (\ref{central}) is  K\"ahler
 iff its  Ricci tensor is 
   $J$-invariant.

\section{Einstein Metrics}

Our first application of the preceding Weyl estimates 
will be to prove new  non-existence 
for Einstein metrics on suitable  smooth compact $4$-manifolds. 
We begin by proving the following technical result: 

\begin{prop}
Let $(M,H)$ be a  polarized smooth compact  
oriented  4-manifold with monopole class $c$. 
Then every $H$-adapted metric satisfies
$$
\frac{1}{4\pi^{2}}\int_{M}\left( \frac{s_{g}^{2}}{24} + 
2|W_{+}|_{g}^{2}\right) d\mu_{g} \geq \frac{2}{3} 
(c_{1}^{+})^{2} .
$$
If equality holds, moreover, $g$ is an almost-K\"ahler metric
of the type described in Proposition \ref{aka}. 
\end{prop}
\begin{proof}
 Let us begin by  rewriting the inequality (\ref{central}) 
 as
 $$\left\|\frac{2}{3}s-2\sqrt{\frac{2}{3}}|W_{+}|\right\|\geq 4\sqrt{2}\pi 
 |c_{1}^{+}|,$$
 where $\|\cdot\|$ denotes the $L^{2}$ norm with respect to $g$.
 By the triangle inequality, we therefore have 
 \begin{equation}
 	\frac{2}{3}\|s\| + {\frac{1}{3}}\| (\sqrt{24}|W_{+}|)\|\geq 4\sqrt{2}\pi 
 |c_{1}^{+}|. 
 	\label{trngl}
 \end{equation}
 We now elect  to interpret the left-hand side as the dot product
 $$(\frac{2}{3} , {\frac{1}{3\sqrt{2}}}) \cdot \left(\|s\|, 
 \|(\sqrt{48}|W_{+}|)\|\right) =  
 \frac{2}{3}\|s\| + {\frac{1}{3}}\| (\sqrt{24}|W_{+}|)\|$$
 in ${\mathbb R}^{2}$. Applying  Cauchy-Schwarz, we thus have 
 $$ \left((\frac{2}{3})^{2} + ({\frac{1}{3\sqrt{2}}})^{2}\right)^{1/2}
 \left(\int_{M}(s^{2}+48|W_{+}|^{2})d\mu \right)^{1/2}
 \geq \frac{2}{3}\|s\| + {\frac{1}{3}}\| (\sqrt{24}|W_{+}|)\|.$$
 Thus
 $$\frac{1}{2}\int_{M}(s^{2}+48|W_{+}|^{2})d\mu \geq 
 32\pi^{2}(c_{1}^{+})^{2},$$
 so that  
 $$\frac{1}{4\pi^{2}}\int_{M}\left( \frac{s_{g}^{2}}{24} + 
2|W_{+}|_{g}^{2}\right) d\mu_{g} \geq \frac{2}{3} 
(c_{1}^{+})^{2} ,$$
as claimed. 
 If equality holds, moreover, (\ref{central}) is saturated, and 
 Proposition \ref{aka} applies. 
\end{proof}

\begin{prop}\label{west}
Let $X$ be a compact oriented $4$-manifold with 
a non-trivial Seiberg-Witten
invariant, and set  
$$M= X\# k \overline{\bcp}_{2} \# \ell (S^{1}\times 
S^{3})$$
for  integers $k, \ell \geq 0$. 
Then any Riemannian metric $g$ on $M$ satisfies
$$
\frac{1}{4\pi^{2}}\int_{M}\left( \frac{s_{g}^{2}}{24} + 
2|W_{+}|_{g}^{2}\right) d\mu_{g} \geq \frac{2}{3} 
(2\chi + 3\tau)(X) ,
$$
and equality is possible only if $k=\ell=0$.  
\end{prop}
\begin{proof}
 The proof is a direct extension of the 
computations in \cite{lebweyl}, although  for $\ell > 0$ we 
now use the holonomy-constrained Seiberg-Witten invariant 
of   \cite{ozsz}; cf. \cite{rio,lebyond}.

Now notice that we may assume that $(2\chi + 3\tau )(X) > 0$,
since otherwise there is nothing to prove. This  has the pleasant 
consequence that any Seiberg-Witten invariant of $X$ is 
independent of polarization, even if $b^{+}(X) =1$. 
Let
 $c_{1}(X)$ denote the first Chern class of a
spin$^{c}$ structure on $X$ for which the Seiberg-Witten invariant
is non-zero, and notice that $(c_{1}(X))^{2}\geq
 (2\chi + 3\tau ) (X) > 0$, because the relevant Seiberg-Witten 
 moduli space must have non-negative virtual dimension.  
 Pull $c_{1}(X)$ back to $M= X\# k \overline{\bcp}_{2}\ell (S^{1}\times 
S^{3})$ via
 the canonical collapsing map, and, by a standard  abuse of 
 notation, let $c_{1}(X)$ also denote this  
 pulled-back class. Thus, with respect to our given polarization,
 $$([c_{1}(X)]^{+})^{2}\geq (c_{1}(X))^{2}\geq
 (2\chi + 3\tau ) (X) > 0.$$
 Now choose  generators
 $E_{1}, \ldots , E_{k}$ for the  pull-backs to $M$ of
  the 
 $k$ relevant copies of $H^{2}(\overline{\bcp}_{2}, {\Bbb Z})$  
 so that 
 $$[c_{1}(X)]^{+}\cdot E_{j}\leq 0, ~~~j = 1, \ldots , k.$$
 Let $\beta_{1}, \ldots, \beta_{\ell}$ be closed curves 
 in $M$ which
 generate of the fundamental groups of  the $\ell$ relevant copies of 
 $S^{1}\times S^{3}$. 
 Then \cite{fs,ozsz} there is a spin$^{c}$ structure on $M$ with
 ${\mathcal S \mathcal W}_{c}(M,[\beta_{1}], \cdots, [\beta_{\ell}]; H)\neq 0$ and 
 $$c_{1}(L)= c_{1}(X) - \sum_{j=1}^{k}  E_{j}.$$
 Thus $c$ is a monopole class of $(M,H)$. 
 But one then has
 \begin{eqnarray*}
 	(c_{1}^{+})^{2} & = & \left([c_{1}(X)]^{+} + 
 	\sum_{j=1}^{k} 	E_{j}^{+}\right)^{2}  \\
 	 & = & ([c_{1}(X)]^{+})^{2}- 2 \sum_{j=1}^{k} [c_{1}(X)]^{+}\cdot  E_{j}
 + (\sum_{j=1}^{k}  E_{j}^{+})^{2}  \\
 	 & \geq  & ([c_{1}(X)]^{+})^{2}  \\
 	 & \geq  & (2\chi + 3\tau ) (X), 
 \end{eqnarray*}
 exactly as claimed.

If equality held, $g$ would be almost-K\"ahler, and $c_{1}(L)$
would be  the anti-canonical class of the associated  almost-complex structure
on $M$. 
But, by construction, 
$$c_{1}^{2}(L)= (c_{1}(X)-\sum E_{j})^{2}= (2\chi + 3\tau ) (X)-k, $$
whereas 
$$(2\chi + 3\tau ) (M) = (2\chi + 3\tau ) (X)-k-4\ell$$
by Mayer-Vietoris. Since $c_{1}^{2}=2\chi + 3\tau$ for
an almost-complex manifold, equality is thus excluded unless 
$\ell =0$. 

There is yet more information available, however. Indeed, 
if equality held, 
the almost-K\"ahler class 
$[\omega ]$ would also necessarily be a non-positive multiple of $c_{1}^{+}$. 
On the other hand, 
our computation shows that equality can only hold if    
 $[c_{1}(X)]^{+}\cdot  E_{j}=0$,
 so it  would follow that   
  $[\omega ]\cdot E_{j}=0$ for all $j$. However, the
 Seiberg-Witten invariant would then also be  non-trivial for a spin$^{c}$ 
 structure with $c_{1}(\tilde{L})= c_{1}(L)+2 E_{1}$,
and a celebrated theorem of Taubes \cite{t2} would then force 
the homology class $E_{j}$ to be represented by a pseudo-holomorphic
$2$-sphere in the symplectic manifold $(M,\omega )$. But the 
(positive!) area of this
sphere with respect to  $g$ would 
then be  exactly $[\omega ]\cdot E_{j}$, contradicting the 
observation 
that $[\omega ]\cdot E_{j}=0$. Thus equality can definitely be excluded
unless $k$ and $\ell$ both vanish. 
\end{proof}

This then gives us a new non-existence theorem for Einstein metrics:

\begin{thm} \label{surg}
Let $X$ be a compact oriented $4$-manifold with 
a non-trivial Seiberg-Witten
invariant. Then 
$$M= X\# k \overline{\bcp}_{2} \# \ell (S^{1}\times 
S^{3})$$
does not admit  Einstein metrics if $k+4\ell \geq 
\frac{1}{3}(2\chi + 3\tau )(X)$. 
\end{thm}
\begin{proof}
We may assume that $(2\chi + 3\tau ) (X) > 0$, since otherwise
the result follows from the Hitchin-Thorpe inequality. 

Now 
$$(2\chi + 3\tau )(M)=\frac{1}{4\pi^{2}}\int_{M}\left( \frac{s_{g}^{2}}{24} + 
2|W_{+}|_{g}^{2} -\frac{|\stackrel{\circ}{r}|^{2}}{2}\right) d\mu_{g}$$
for any metric on $g$ on $M$.  If $g$ is an Einstein metric,
the trace-free part $\stackrel{\circ}{r}$ of the Ricci curvature
vanishes, and  we then have
\begin{eqnarray*}
(2\chi + 3\tau )(X) -k-4\ell 	 & = & (2\chi + 3\tau )(M)  \\
	 & = &  \frac{1}{4\pi^{2}}\int_{M}\left( \frac{s_{g}^{2}}{24} + 
2|W_{+}|_{g}^{2}\right) d\mu_{g} \\
	 & \geq & \frac{2}{3} (2\chi + 3\tau )(X)  \\
\end{eqnarray*}
 by Proposition \ref{west}, with  equality only if $k$ and $\ell$
 both vanish.  If $M$ carries an Einstein metric, it therefore follows 
 that 
 $$ \frac{1}{3}(2\chi + 3\tau )(X) > k+ 4\ell .$$
 The claim now follows by contraposition. 
\end{proof}

Specializing to the case of complex surfaces, we  now have: 

\begin{cor}\label{blowup}
Let $(X,J_{X})$ be a minimal complex surface of general type, and 
let $(M,J_{M})$ be obtained from $X$ by blowing up 
$k$  points. Then the smooth compact $4$-manifold $M$ 
does not admit any Einstein metrics if 
 $k \geq 
\frac{1}{3}c_{1}^{2}(X)$. 
\end{cor}

\begin{xpl} Let $X\subset {\bcp}_{3}$ be a hypersurface of 
degree $6$. Since the canonical class on $X$ is twice the
hyperplane class, $c_{1}^{2}(X)= 2^{2}\cdot 6= 24$.  Corollary 
\ref{blowup}
therefore tells us that if we  
blow up $X$  at $8$ points, the  resulting $4$-manifold
$$M = X\# 8 \overline{\bcp}_{2}$$
does not admit Einstein metrics. 

Now let us compare the complex surface $M$
with the {\em Horikawa surface} $N$ 
obtained as a ramified double  cover of 
${\mathbb C \mathbb P}_{1}\times {\mathbb C \mathbb P}_{1}$ 
 branched at a generic curve of bidegree 
 $(6,12)$. Both of these simply connected complex surfaces
 have
 $c_{1}^{2}=16$ and $p_{g}=10$, and 
 the underlying oriented $4$-manifolds therefore have
 $b_{+}= 21$ and $b_{-}=93$. In particular, 
 both have signature $\tau = -72 \not \equiv 0 \bmod 16$,
 so, by Rochlin's theorem, neither is spin. 
 Thus $M$ and $N$ have
 isomorphic intersection
 forms  by the Minkowski-Hasse classification, and
 are therefore homeomorphic by Freedman's theorem \cite{freedman}. 
However, $N$ has ample canonical line bundle, and 
so admits a K\"ahler-Einstein metric by Yau's theorem. 
Thus, although $M$ and $N$ are homeomorphic, one admits
Einstein metrics, while the other doesn't.

If we instead start with a hypersurface  
$X\subset {\bcp}_{3}$ of degree $\geq 7$, one 
can construct homeomorphic pairs $(M,N)$ for
which the minimal complex surface sits well above the
Noether line. Details are left to the interested reader.   
\end{xpl}

\begin{rmk}
Corollary \ref{blowup} is the direct descendant  of 
an analogous result in \cite{lno}, where, using only
scalar curvature estimates, a similar conclusion was
proved for $k \geq \frac{2}{3} c_{1}^{2}(X)$. It 
was later pointed out by Kotschick \cite{kot} that such a result alone
suffices to imply the existence of homeomorphic pairs
consisting of an Einstein manifold and a 
$4$-manifold which does not admit
 Einstein metrics;  
 however,  the   examples that arise by this method
are quite complicated, and  have huge $c_{1}^{2}$.   
 
 The intermediate step between \cite{lno} and 
 Corollary \ref{blowup} may be found in \cite{lebweyl},
 where  Seiberg-Witten estimates of Weyl curvature
  were first introduced. While crude by present standards,
 the method used there 
  did lead to an obstruction when $k \geq \frac{25}{57} c_{1}^{2}(X)$,
  or about two-thirds of the way to the present result. 
  
  The fact that $S^{1}\times S^{3}$ handles can be added
  to a $4$-manifold 
  without losing Seiberg-Witten control of the scalar curvature was 
  first observed by Petean \cite{jp}, although Seiberg-Witten theory
  only enters his result in an indirect manner.
  The search for a 
   direct Seiberg-Witten explanation of this phenomenon 
   then led the   author and  del Rio \cite{lebyond,rio}
   to a discovery of the generalized  Seiberg-Witten invariant
   used here. 
 In fact, however, this turned out to    merely be
  a {\em re}-discovery;   
   quite   different considerations had already
  led   Ozsv\'ath and ~Szab\'o \cite{ozsz}
   to  develop a full-blown theory, described in a
   preprint written the previous  month.
\end{rmk}

Let us now recall that the Hitchin-Thorpe inequality
\cite{hit,thorpe, bes} asserts that 
if a compact
oriented $4$-manifold $M$ 
admits an Einstein metric, then  
$$\frac{|\tau (M)|}{\chi(M)}\leq \frac{2}{3}.$$
Hitchin went on to observe that the converse 
is certainly false, because the argument
 shows in particular that a simply connected 
 $4$-manifold with 
$2\chi(M) = -3\tau (M)$
can admit an Einstein metric only if its a 
 $K3$ surface. Examples of $4$-manifolds with
$|\tau (M)|/\chi(M) < 2/3$ which do not admit   
Einstein metrics were first
constructed by Gromov \cite{grom}, using 
his simplicial 
 volume invariant; however, this method only works in   
 the presence of  an infinite  fundamental group. 
In \cite{lno},  simply connected examples were  first 
constructed, using 
Seiberg-Witten estimates for the scalar curvature. 
Shortly thereafter, 
Sambusetti \cite{samba} showed that 
the entropy estimates of Besson-Courtois-Gallot \cite{bcg2}
allow one to construct examples of arbitrary
Euler characteristic and signature, but again with 
huge fundamental group; see Petean \cite{jp} and 
del Rio \cite{rio} for similar results.  
In the {\em simply connected} case, however, 
it remains to be seen whether
$(\chi , \tau)$ can really be arbitrarily specified,
subject to the obvious constraints  $\chi \geq 2+
|\tau |$ and $\chi \equiv \tau \bmod 2$. Indeed, even 
the following is open:

\begin{question}
For every $q\in (-1,1)\cap {\mathbb Q}$, are there  smooth, 
compact  simply connected $4$-manifolds with 
$\tau /\chi  =q$ 
which do not admit 
Einstein metrics? 
\end{question}

Hitchin-Thorpe  gives  one a resoundingly affirmative answer in the range  
$\frac{2}{3}\leq |q| < 1$. The present results 
allow us to  improve this as follows:  

\begin{cor}
Let $q$ be a rational number with $\frac{8}{23}\leq | q |< 1$. 
Then there 
are smooth, compact,  simply connected $4$-manifolds  with 
${\tau}/{\chi}=q$
which do not admit 
Einstein metrics.
\end{cor} 
\begin{proof}
For $m$  any even integer bigger than $17$ million,
Chen \cite{chen} has constructed a 
simply connected minimal 
complex surface $X$ of general type 
(in fact, a hyperelliptic fibration) 
with $\tau (X)=m$ and $\chi(X)=4m$.
If we now blow up such a surface at 
$k$ points, where $k\geq \frac{11}{3}m$,
the resulting simply connected $4$-manifold $M = X\# k \overline{\bcp}_{2}$
does not admit Einstein metrics by Corollary \ref{blowup}. 
Now for this manifold we have  
$$\frac{\tau (M)}{\chi(M)}= \frac{5}{4+\frac{k}{m}}-1,$$
and the right-hand side sweeps out  ${\mathbb Q}\cap  
(-1, -\frac{8}{23}]$ as $\frac{k}{m}$ ranges over
${\mathbb Q}\cap [ \frac{11}{3}, \infty)$. 
  This proves the claim for $q$ negative. 
The case of $q$ positive then follows by reversing  
orientation. 
\end{proof}

The Weyl estimates of \S \ref{weylest}  also have
interesting ramifications for the theory of {\em anti-self-dual}
$4$-manifolds. Recall that an oriented
 Riemannian $4$-manifold is said to be
anti-self-dual if it satisfies the $W_{+}\equiv 0$, and that this
condition is conformally invariant. Compact anti-self-dual 
manifolds exist in profusion. Indeed,  Taubes \cite{tasd} has shown 
that for any smooth compact orientable $X^{4}$, there is an integer
$k_{0}$ such that $M=X\# k\overline{\bcp}_{2}$ admits metrics
with $W_{+}=0$ provided that $k \geq k_{0}$. In particular, 
if we blow up a symplectic $4$-manifold with $b^{+}> 1$ 
 at enough
points, we obtain a manifold which  admits
anti-self-dual metrics and also has a non-trivial Seiberg-Witten 
invariant. Let us now consider the scalar curvature 
of such manifolds.

\begin{lem}\label{asdm}
Let $(M,g)$ be a compact {\em anti-self-dual} 
4-manifold with a non-zero Seiberg-Witten invariant.  Then 
 $$\int_{M}s^{2}d\mu \geq 
 72\pi^{2}(c_{1}^{+})^{2},$$
 where $c_{1}^{+}$ is again the self-dual part of the 
 first Chern class of the relevant spin$^{c}$ structure.
 If equality holds, moreover, $g$ is an almost-K\"ahler metric
 of the type describe in  Proposition \ref{aka}.
 \end{lem}
 \begin{proof}
 If $W_{+}\equiv 0$, 
 (\ref{central}) becomes 
$$	\int_{M}\frac{4}{9} s^{2} 
	d\mu 
\geq 32\pi^{2}(c_{1}^{+})^{2},
$$ 
and the desired inequality is therefore an immediate consequence. 
Moreover, if equality holds, Proposition \ref{aka} applies,
and  $g$ is an almost-K\"ahler manifold with the relevant 
special properties.
 \end{proof}
 
 Notice that this result may be used as a vanishing theorem. For example, 
 if $c$ is a spin$^{c}$ structure on a hyperbolic $4$-manifold,
 the corresponding Seiberg-Witten invariant will vanish 
 unless $(c_{1}^{+})^{2} \leq \frac{8}{3}\chi$. Since
 the non-negativity of the virtual dimension of 
 the Seiberg-Witten moduli space  would require that
   $(c_{1}^{+})^{2} \geq 2\chi + |c_{1}^{-}|^{2}$,
 this would seem to lend some weak support to the conjecture that 
 the Seiberg-Witten
 invariants of a 
 hyperbolic $4$-manifold are all zero. 
 In any case, it also implies the following result:
 
 \begin{thm}
Let $X$ be a compact oriented $4$-manifold with 
a non-trivial Seiberg-Witten
invariant. Then 
$$M= X\# k \overline{\bcp}_{2} \# \ell (S^{1}\times 
S^{3})$$
does not admit anti-self-dual Einstein metrics if $k+4\ell \geq 
\frac{1}{4}(2\chi + 3\tau )(X)$. 
 \end{thm}
 \begin{proof}
 As in the proof of Proposition \ref{west},  
 there is a monopole class $c$ of $(M,H)$ with
 $$c_{1}(L)= c_{1}(X) - \sum_{j=1}^{k}  E_{j},$$
and  such that 
 $$
 (c_{1}^{+})^{2} \geq (2\chi + 3\tau ) (X).
 $$
 If $M$ carried an anti-self-dual Einstein metric,  Lemma \ref{asdm} 
 would 
 then tell us that
 \begin{eqnarray*}
 (2\chi + 3\tau )(X) -k-4\ell	 & = & (2\chi + 3\tau )(M)
  \\
 	 & = & \frac{1}{4\pi^{2}}\int_{M}\frac{s^{2}}{24}d\mu
 	  \\
 	 & \geq  & \frac{72\pi^{2}}{96\pi^{2}}(c_{1}^{+})^{2}
 	  \\
 	 & \geq  & \frac{3}{4}(2\chi +3\tau )(X), 
 \end{eqnarray*}
 so that 
 $$k+4\ell \leq \frac{1}{4}(2\chi +3\tau )(X).$$
 Moreover, equality could only hold if 
 $k$ and $\ell$ both vanished for precisely the same reasons 
 delineated in the proof of Theorem \ref{surg}. 
 
The result now follows by contraposition. 
\end{proof}

\section{Ricci Curvature  Estimates}

One of the most striking aspects of Seiberg-Witten theory is that
it allows one to calculate the precise value of the invariant 
$${\mathcal I}_{s}(M)=\inf_{g}\int_{M}s^{2}_{g}d\mu_{g}$$
for large classes of $4$-manifolds $M$; for example, if 
$M$ is a complex surface of general type, and if
$X$ is its minimal model, then \cite{lno} 
\begin{equation}
	\inf_{g}\int_{M}s^{2}_{g}d\mu_{g} = 32\pi^{2}c_{1}^{2}(X).
	\label{scal}
\end{equation}
However, one would also like to compute the 
corresponding infimum 
\begin{equation}
{\mathcal I}_{r}(M)=	\inf_{g} \int_{M}|r_{g}|^{2}d\mu_{g}
	\label{rich}
\end{equation}
for the {Ricci} curvature 
$r$. Since 
$$\int |r|^{2}d\mu = \int \left(\frac{s^{2}}{4}+ 
|\stackrel{\circ}{r}|^{2}\right)d\mu ,$$
where $\stackrel{\circ}{r}$ is the trace-free part of the
Ricci curvature, we obviously have
\begin{equation}
	{\mathcal I}_{r}(M) \geq \frac{1}{4}{\mathcal I}_{s}(M)
	\label{taut}
\end{equation}
and the scalar curvature estimate obviously 
gives us some important information about this problem.
Moreover, because K\"ahler-Einstein metrics saturate the 
lower bound for the scalar curvature, one actually gets
equality in (\ref{taut}) if $M$ happens to be a 
{\em minimal} complex surface of general type.   
For non-minimal surfaces, however, the situation is
dramatically different; indeed, even the relatively crude Weyl estimates
of \cite{lebweyl} allow one to deduce that, for any 
complex surface $M$ of general type,  
\begin{center}{\em 
equality holds in (\ref{taut}) $\Longleftrightarrow$ $M$ is minimal.} 
\end{center} 
Until now, however, non-trivial direct calculations of ${\mathcal 
I}_{r}$ 
have remained beyond the scope of the theory. Fortunately, 
we will now see that the estimates of \S \ref{weylest}
are perfectly suited to this purpose. The key observation is the
following:

\begin{prop} \label{glove}
Let $(M,H)$ be a  polarized smooth compact  
oriented  4-manifold with monopole class $c$. 
Then every $H$-adapted metric satisfies
$$
\frac{1}{4\pi^{2}}\int_{M}\left( \frac{s_{g}^{2}}{24} + 
\frac{1}{2}|W_{+}|_{g}^{2}\right) d\mu_{g} \geq \frac{1}{2} 
(c_{1}^{+})^{2} .
$$
Moreover, equality holds iff $g$ is a K\"ahler
metric of constant scalar curvature, and the 
spin$^{c}$ structure is the one determined by 
a $g$-compatible complex structure on $M$. 
\end{prop}
\begin{proof}
We again begin with  inequality (\ref{trngl}),
 $$
 	\frac{2}{3}\|s\| + {\frac{1}{3}}\| (\sqrt{24}|W_{+}|)\|\geq 4\sqrt{2}\pi 
 |c_{1}^{+}|,
 $$
but this time interpret the left-hand side as the dot product
 $$(\frac{2}{3} , {\frac{\sqrt{2}}{3}}) \cdot \left(\|s\|, 
 \|(\sqrt{12}|W_{+}|)\|\right) $$
 in ${\mathbb R}^{2}$. Applying the Cauchy-Schwarz inequality, we 
 then 
 obtain 
 $$ \left((\frac{2}{3})^{2} + ({\frac{\sqrt{2}}{3}})^{2}\right)^{1/2}
 \left(\int_{M}(s^{2}+12|W_{+}|^{2})d\mu \right)^{1/2}
 \geq \frac{2}{3}\|s\| + {\frac{1}{3}}\| (\sqrt{24}|W_{+}|)\|.$$
 Thus
 $$\frac{2}{3}\int_{M}(s^{2}+12|W_{+}|^{2})d\mu \geq 
 32\pi^{2}(c_{1}^{+})^{2},$$
 and 
 $$\frac{1}{4\pi^{2}}\int_{M}\left( \frac{s_{g}^{2}}{24} + 
\frac{1}{2}|W_{+}|_{g}^{2}\right) d\mu_{g} \geq \frac{1}{2} 
(c_{1}^{+})^{2} ,$$
as claimed.

 If equality holds, our use of  Cauchy-Schwarz forces 
 $$(\frac{2}{3}:{\frac{\sqrt{2}}{3}})=(\|s\|:\|(\sqrt{12}|W_{+}|)\|),$$
which is to say that 
 $$\int |W_{+}|^{2}d\mu = \int \frac{s^{2}}{24}d\mu.$$
 In this case, we  then have 
 $$\int_{M} s_{g}^{2} d\mu_{g} = 32 \pi^{2}
(c_{1}^{+})^{2} ,$$
and it then follows \cite{lpm} that the metric is  K\"ahler,
and has constant scalar curvature.  Conversely, a K\"ahler
surface with constant scalar curvature satisfies $s\equiv \sqrt{24}|W_{+}|$
and $\int s^{2}d\mu =  32 \pi^{2}
(c_{1}^{+})^{2}$, so that equality is attained for any such 
metric.
\end{proof}

 This then gives rise to a Ricci curvature estimate
 via  Gauss-Bonnet  formul{\ae} : 

\begin{prop}
Let $(M,H)$ be a  polarized smooth compact  
oriented  4-manifold with monopole class  $c$. 
Then every $H$-adapted metric satisfies
$$
\frac{1}{8\pi^{2}}\int_{M}|r_{g}|^{2} d\mu_{g} \geq  
2(c_{1}^{+})^{2} - (2\chi + 3\tau )(M),
$$
where $r$ denotes the Ricci curvature of $g$. 
\end{prop}
\begin{proof}
Since (\ref{bonnet}) tells us that 
$$(2\chi + 3\tau )(M)=\frac{1}{4\pi^{2}}\int_{M}\left( \frac{s_{g}^{2}}{24} + 
2|W_{+}|_{g}^{2} -\frac{|\stackrel{\circ}{r}|^{2}}{2}\right) 
d\mu_{g} , $$ 
it follows that 
\begin{equation}
	\frac{1}{8\pi^{2}}\int_{M} |r|^{2}d\mu = \frac{1}{\pi^{2}}
\int_{M}\left( \frac{s_{g}^{2}}{24} + 
\frac{1}{2}|W_{+}|_{g}^{2}\right)d\mu_{g} - (2\chi + 3\tau )(M)
	\label{gbs}
\end{equation}
for any metric on $g$ on $M$. The result now follows immediately from 
Proposition \ref{glove}. 
\end{proof}

This observation  supplies the key to  the main result
of this section: 

\begin{thm} \label{l2ric}
Let $X$ be any minimal complex surface of general type,
and let $M= X\# k \overline{\mathbb C \mathbb P}_{2} \# \ell (S^{1}\times 
S^{3})$. Then 
$${\mathcal I}_{r}(M)=
 \inf_{g} \int_{M} |r|^{2}d\mu = 8\pi^{2} (c_{1}^{2}(X) + k+4\ell) .$$
\end{thm}
\begin{proof}
Since 
$$(2\chi + 3\tau )(M)= c_{1}^{2}(X)-k-4\ell,$$
equation (\ref{gbs}) tells us that it suffices to prove that
$$
\inf_{g} \frac{1}{4\pi^{2}}
\int_{M}\left( \frac{s_{g}^{2}}{24} + 
\frac{1}{2}|W_{+}|_{g}^{2}\right)d\mu_{g} = \frac{1}{2}c_{1}^{2}(X).
$$

Now  observe that the proof of 
Proposition \ref{west} tells us that that, for every
metric $g$ on $M$, there is a spin$^{c}$ structure $c$
with non-zero Seiberg-Witten invariant such that 
$$(c_{1}^{+})^{2} \geq  (2\chi + 3\tau ) (X) = c_{1}^{2}(X).$$
Hence 
$$
\frac{1}{4\pi^{2}}
\int_{M}\left( \frac{s_{g}^{2}}{24} + 
\frac{1}{2}|W_{+}|_{g}^{2}\right)d\mu_{g} \geq \frac{1}{2}c_{1}^{2}(X)$$
for every metric $g$ on $M$.
It follows  that 
$$
\inf_{g} \frac{1}{4\pi^{2}}
\int_{M}\left( \frac{s_{g}^{2}}{24} + 
\frac{1}{2}|W_{+}|_{g}^{2}\right)d\mu_{g} \geq \frac{1}{2}c_{1}^{2}(X).
$$

To finish the proof, it  suffices to produce a family 
of metrics $g_{\varepsilon}$ on $M$ such that 
$$
\limsup_{\varepsilon \to 0^{+}}  \frac{1}{4\pi^{2}}
\int_{M}\left( \frac{s_{g_{\varepsilon}}^{2}}{24} + 
\frac{1}{2}|W_{+}|_{g_{\varepsilon}}^{2}\right)d\mu_{g_{\varepsilon}} 
\leq
\frac{1}{2}c_{1}^{2}(X).$$
To do this, let $\check{X}$ be the pluri-canonical model 
of $X$, which carries a K\"ahler-Einstein orbifold metric
$\check{g}$ by an immediate generalization \cite{rkob,tsuj} 
of the Aubin/Yau
solution  \cite{aubin,yau}  of the $c_{1}<0$ case of
the Calabi conjecture.  Now
$$\int_{\check{X}}\left( \frac{s_{\check{g}}^{2}}{24} + 
\frac{1}{2}|W_{+}|_{\check{g}}^{2}\right)d\mu_{\check{g}}
= \frac{1}{2}c_{1}^{2}(\check{X})= \frac{1}{2}c_{1}^{2}(X)$$
because $|W_{+}|^{2}\equiv s^{2}/24$ for any K\"ahler 
surface, whereas a K\"ahler-Einstein surface has 
Ricci form  given by $s \omega/4$, where $\omega$ is the 
K\"ahler form. Let $o_{1}, \ldots , o_{n}$ be the orbifold points
of $\check{X}$, if there are any. 
Each of these singularities is of type A-D-E, meaning that 
$o_{j}$ has a neighborhood   is modeled on ${\mathbb C}^{2}/\Gamma_{j}$ for
some discrete sub-group $\Gamma_j$ of $SU(2)$. 
The passage from $\check{X}$ to $X$ is accomplished by 
replacing each $o_{j}$ by a collection of $(-2)$-curves
with intersections
 determined by the Dynkin diagram (Coxeter graph)
of type A-D-E corresponding to $\Gamma_{j}$. 

One can now construct smooth metrics on the minimal model
$X$ by modifying the orbifold metric $\check{g}$, without introducing
substantial amounts of extra volume, Ricci curvature, or 
self-dual Weyl curvature. 
Indeed, for each $\Gamma_{j}$ there are \cite{kron} gravitational 
instanton metrics on  the minimal resolution of
${\mathbb C}^{2}/\Gamma_{j}$, which is precisely obtained by replacing the
origin with a set of  $(-2)$-curves as above.
These gravitational instanton metrics are Ricci-flat and 
anti-self-dual, 
and they closely approximate the 
  Euclidean metric on ${\mathbb C}^{2}/\Gamma_{j}$, 
  being of the form $\delta_{jk}+O(\varrho^{-4})$
  in the asymptotic region,
  where 
 $\varrho$ is the Euclidean radius; moreover, the first and second
  derivatives of such a metric are of order 
   $O(\varrho^{-5})$ and  $O(\varrho^{-6})$,
  respectively, in these coordinates. Fix such a metric, and consider the 
  family of gravitational instanton metrics $h_{\varepsilon}$ obtained by 
  multiplying it by $\varepsilon^{4}$, and then making the homothetic
 change $x\mapsto \varepsilon^{2}x$ of asymptotic coordinates.
 The metric $h_{\varepsilon}$ is then uniformly 
   $\delta_{jk} + O(\varepsilon^{8} \varrho^{-4})$ 
   on the complement of the ball of radius $\varepsilon^{2}$,
   with first and second derivatives
   $O(\varepsilon^{8} \varrho^{-5})$
   and   $O(\varepsilon^{8} \varrho^{-6})$,
   respectively. 
   Take geodesic spray coordinates about each 
   of the orbifold points $o_{j}$, so that  
   $\check{g}$ takes the form $\delta_{jk}+ O(\varrho^{2})$
   in these coordinates. 
   Delete the $\varepsilon$ ball around each 
   orbifold point, and replace it in the $\varrho \leq 
   \varepsilon$ region of the corresponding gravitational instanton;
   and on the  transitional annulus $\varrho \in [\varepsilon ,2 \varepsilon ]$,
   consider the metric
   \begin{equation}
   	g_{\varepsilon}= 
 f (\frac{\varrho}{\varepsilon}) h_{\varepsilon}+
 [1- f (\frac{\varrho}{\varepsilon})]
  \check{g} , 
   	\label{trans}
   \end{equation}
  where  $f: [1,2] \to [0,1]$ is 
  a fixed smooth function 
  which is $\equiv 1$ near $1$ and $\equiv 0$ near $2$.
  This gives us a family of metrics on $X$, the geometry of which 
  will be analyzed in a moment.  

However, what we really need is a metric
on $M$ rather than on $X$, so we need to 
 modify $\check{g}$ a bit more. 
Now $S^{1}\times S^{3}$ and
$\overline{\mathbb C \mathbb P}_{2}$ both admit anti-self-dual metrics
of positive scalar curvature --- namely, the obvious product 
 metric on  $S^{1}\times S^{3}$, and the 
Fubini-Study metric on the reverse-oriented complex
projective plane $\overline{\mathbb C \mathbb P}_{2}$. 
The complement of any point in 
$S^{1}\times S^{3}$ or 
$\overline{\mathbb C \mathbb P}_{2}$ therefore  admits 
an asymptotically flat metric with $s\equiv 0$ and $W_{+}\equiv 0$,
namely the standard metric rescaled by the Yamabe Greens function 
of the appropriate point. One can then find asymptotic
charts for these manifolds, parameterized by the complement of, say,  
the unit  
 ball in  ${\mathbb R}^{4}$, in which these metrics
 differ from the Euclidean metric by $O(\varrho^{-2})$ and 
 in which first and second partial derivatives are respectively
 $O(\varrho^{-3})$ and $O(\varrho^{-4})$. 
 (Notice the fall-off  rate is slower than in the previous case;
 these metrics have positive mass!)  Again multiplying 
 by $\varepsilon^{4}$, where $\varepsilon$ is to be viewed as
 an auxiliary parameter,  and making a homothety as before, 
 we thus get metrics
 on the complement of the $\varepsilon^{2}$ ball in  ${\mathbb R}^{4}$ 
 of the form $\delta_{jk}+ O(\varepsilon^{4}\varrho^{-2})$, with 
 first and second partial derivatives 
 of order $O(\varepsilon^{4}\varrho^{-3})$
 and $O(\varepsilon^{4}\varrho^{-4})$, respectively. 
 Again, let $h_{\varepsilon}$ denote the restriction of the resulting 
 metric on the the $ \varrho > 2\varepsilon$ regions of $k$ disjoint punctured 
 copies of $\overline{\mathbb C \mathbb P}_{2}$
 and $\ell$ disjoint punctured copies of $S^{1}\times S^{3}$.
 Take
 geodesic spray coordinates around non-singular 
 points $p_{1}, \ldots, p_{k}, q_{1}, \ldots , q_{\ell}$ of
 $\check{X}$, and delete $\varepsilon$-balls around each one, and 
  define a metric $g_{\varepsilon}$ on $M$ as
 $h_{\varepsilon}$ when $\varrho < \varepsilon$, 
 as $\check{g}$ when $\varrho > 2\varepsilon$, and as
 (\ref{trans}) on the transitional annuli 
 $\varrho\in [\varepsilon , 2 \varepsilon]$. 
 The sectional curvature of the metrics $g_{\varepsilon}$
  are then uniformly  bounded   on the annuli
  $\varepsilon \leq \varrho \leq  
2 \varepsilon $ as  $\varepsilon\searrow 0$, while
 the volumes of these annuli are
$O(\varepsilon^{4})$ as $\varepsilon \to 0$.
Since the metrics we have used to replace the $\varepsilon$-balls
are all scalar-flat and anti-self-dual, we therefore have 
$$
\frac{1}{4\pi^{2}}
\int_{M}\left( \frac{s_{g_{\varepsilon}}^{2}}{24} + 
\frac{1}{2}|W_{+}|_{g_{\varepsilon}}^{2}\right)d\mu_{g_{\varepsilon}}
\leq  C\varepsilon^{4} +\frac{1}{4\pi^{2}}\int_{\check{X}}\left( 
\frac{s_{\check{g}}^{2}}{24} + 
\frac{1}{2}|W_{+}|_{\check{g}}^{2}\right)d\mu_{\check{g}} ,
$$
and it thus follows that 
$$
\limsup_{\varepsilon \to 0^{+}}  \frac{1}{4\pi^{2}}
\int_{M}\left( \frac{s_{g_{\varepsilon}}^{2}}{24} + 
\frac{1}{2}|W_{+}|_{g_{\varepsilon}}^{2}\right)d\mu_{g_{\varepsilon}} 
\leq 
\frac{1}{2}c_{1}^{2}(X).$$
The theorem follows.
\end{proof}

\begin{cor}
Let $M$ be the underlying $4$-manifold of any compact complex surface of 
Kodaira dimension $\geq 0$. Let $X$ be the minimal model of $M$, 
and let $k$ be the number of points at which $X$ must be blown up 
so as  to obtain $M$. Then 
$$ \inf_{g} \int_{M} |r|^{2}d\mu = 8\pi^{2} (c_{1}^{2}(X) + k) 
,$$
where the infimum is taken over all smooth Riemannian metrics on $g$ on $M$. 
\end{cor}
\begin{proof}
Complex surfaces of  Kodaira dimension  $2$, are, 
by definition, precisely those of general type; for these, the result 
is just the $\ell=0$ case of Theorem \ref{l2ric}. On the other hand,
the analogous assertion for Kodaira dimensions $0$ and $1$ immediately 
follows
 from the fact \cite{lky} that 
 any elliptic surface admits sequences of metrics for which
  $\mv (M,g)\searrow 0$, but for which $s$ and $|W_{+}|$ are uniformly  bounded. 
\end{proof}

In particular, in view of (\ref{scal}), 
this gives a satisfying new proof of the following:

\begin{cor}
Let $M$ be the underlying $4$-manifold of any
 complex surface of Kodaira dimension $\geq 0$. Then
\begin{center}
equality holds in (\ref{taut}) $~\Longleftrightarrow ~$
 $M$ is minimal. 
\end{center}
\end{cor}

We are also see that the invariant ${\mathcal I}_{r}$
is extremely sensitive to changes in differentiable structure
on a topological $4$-manifold: 
\begin{cor}
Let $X$ and $\tilde{X}$ be two minimal simply connected complex 
surfaces with the same geometric genus $p_{g}\neq 0$,
but with  $c_{1}^{2}(X) - c_{1}^{2}(\tilde{X})= j > 0$. 
Then for all positive integers $k$, the 
$4$-manifolds $\tilde{X}\# k\overline{\bcp}_{2}$ and 
$X\# (j+k)\overline{\bcp}_{2}$ are homeomorphic,
but have different invariants ${\mathcal I}_{r}$.
\end{cor}
\begin{proof}
The simply connected $4$-manifolds
$\tilde{X}\# k\overline{\bcp}_{2}$ and 
$X\# (j+k)\overline{\bcp}_{2}$ are 
non-spin and have the same invariants 
$b_{\pm}$; Freedman's classification \cite{freedman} 
therefore tells us that they are homeomorphic. 
On the other hand, because $p_{g}\neq 0$,
these complex surfaces  have Kodaira dimension 
$\geq 0$. Theorem \ref{l2ric} therefore tells us that 
$$\frac{1}{8\pi^{2}}
{\mathcal I}_{r}(\tilde{X}\# k\overline{\bcp}_{2})=
c_{1}^{2}(\tilde{X})+k= c_{1}^{2}(X)-j+k$$
whereas 
$$
\frac{1}{8\pi^{2}}
{\mathcal I}_{r}({X}\# (j+k)\overline{\bcp}_{2})=
 c_{1}^{2}(X)+(j+k). 
$$
Thus these manifolds have unequal invariants 
${\mathcal I}_{r}$, as claimed. 
\end{proof}

Note that pairs $(X,\tilde{X})$ of the above kind 
are as common as garden weeds; cf. e.g. \cite{bpv}. For  example, if $X$
is of general type, we can always find 
a simply connected properly elliptic surface $\tilde{X}$ with the
same geometric genus as $X$.

\section{Sectional Curvature and Volume}

In the previous section, we considered  a differential-topological invariant 
arising out of the  $L^{2}$-norm of 
Ricci curvature.  We  will now
turn to a related problem which arises in connection with the  
 $L^{\infty}$ norm of sectional curvature.

 Recall that the sectional curvature $K(P)$ of 
 a Riemannian $n$-manifold $(M,g)$ is 
  a smooth function on the Grassmann bundle 
 of $2$-dimensional subspaces $P$ of $TM$.
 At each point  $x\in M$,   let us consider the  {\em bottom 
  sectional curvature}, given by    
$$ \underline{K}(x):= \min_{P\subset T_{x}M} K(P).$$
Then   $\underline{K}: M\to {\Bbb R}$ is automatically
a Lipschitz continuous  function, 
 although it need not be differentiable in general. 
  Notice that
 we tautologically have $$\underline{K}\leq \frac{s}{n(n-1)},$$ so 
 that  
  $\underline{K}$ is  negative
 on  $\{ x\in M~|~ s(x) < 0\}$.

 \begin{lem}\label{wgl}
 Let $(M,g)$ be an oriented Riemannian $4$-manifold,
 and, for $x$ in $M$, let $w(x)\leq 0$ denote the smallest
 eigenvalue of $W_{+}(x): \Lambda^{+}_{x}\to \Lambda^{+}_{x}$.
 Then 
 $$\underline{K}(x) \leq \frac{s(x)}{12}+ \frac{w(x)}{2}.
 $$
 If equality holds at $x$, then $W_{-}(x)=0$;
 and the converse also holds if $g$ is Einstein. 
 \end{lem}
\begin{proof}
First recall that, in terms of the decomposition
$$\Lambda^{2}= \Lambda^{+}\oplus \Lambda^{-},$$
the curvature operator of an oriented Riemannian 
$4$-manifold can be expressed in the form 
  \begin{equation}
\label{curv}
{\cal R}=
\left(
\mbox{
\begin{tabular}{c|c}
&\\
$W_++\frac{s}{12}$&$B$\\ &\\
\cline{1-2}&\\
$B^{*}$ & $W_{-}+\frac{s}{12}$\\&\\
\end{tabular}
} \right) , 
\end{equation}
where $W_{+}$ and $W_{-}$ are trace-free, and 
where $B: \Lambda^{-}\to \Lambda^{+}$ exactly corresponds
to $\stackrel{\circ}{r}$.  
On the other hand, 
every 2-form $\psi$ on $M$ can be uniquely written as
$\psi = \psi^{+} + \psi^{-}$,
where $\psi^{\pm}\in \Lambda^{\pm}$. Now a 
2-form is expressible as a simple wedge product of 
1-forms iff $\psi \wedge \psi =0$.
Thus $\psi$ is simple and has unit length iff
$|\psi^{+}|^{2}=|\psi^{-}|^{2}=\frac{1}{2}$.
The sectional curvature in the corresponding $2$-plane is then
\begin{equation}
	K(P)= \frac{s}{12}+ \langle \psi^{+} , W_{+}(\psi^{+})\rangle  +\langle 
\psi^{-}, W_{-}(\psi^{-})\rangle + 2 \langle 
\psi^{+}, B\psi^{-}\rangle .
	\label{sec}
\end{equation}
Notice, however, that the last term changes sign if 
$\psi^{-}$ is replaced by $-\psi^{-}$; geometrically,
this means that the average of the sectional curvatures of an
orthogonal pair of $2$-planes only depends on the scalar and
Weyl curvatures. Thus,  there is always a $2$-plane  $P^{\prime}$ for which 
the sectional curvature is exactly
$$K(P^{\prime}) =  \frac{s}{12}+ \langle \psi^{+} , 
W_{+}(\psi^{+})\rangle  +\langle 
\psi^{-}, W_{-}(\psi^{-})\rangle .$$
Hence
\begin{eqnarray}
	\underline{K} & \leq  & \frac{s}{12}+ 
	\inf_{|\psi^{+}|^{2}=\frac{1}{2}}\langle \psi^{+} , 
W_{+}(\psi^{+})\rangle  + \inf_{|\psi^{-}|^{2}
=\frac{1}{2}}\langle 
\psi^{-}, W_{-}(\psi^{-})\rangle
	\nonumber  \\
	 & = &   \frac{s}{12}+ 
\frac{w}{2}+\frac{\tilde{w}}{2},
	\label{crafty}
\end{eqnarray}
where $\tilde{w}\leq 0$ is the lowest eigenvalue of $W_{-}$. 
In particular, 
$$\underline{K}\leq \frac{s}{12}+ 
\frac{w}{2},$$
as claimed.
\end{proof}

 \begin{prop}\label{lower} 
 Let $(M,g)$ be a compact oriented Riemannian $4$-manifold
  on which there is a solution of the Seiberg-Witten 
  equations. Let 
 $c_{1}(L)$ be the  first Chern class of the relevant spin$^{c}$ structure,
 and let $c_{1}^{+}$ denote its self-dual part with respect to $g$.
  Then
 $$-V^{1/2}\min_{x\in M} \left( \underline{K}(x)+\frac{s(x)}{12}\right)  
 \geq {\sqrt{2}} {\pi}|c_{1}^{+}|  ,
 $$
 where  $V=\mv (M,g)$ denotes the total volume of $(M,g)$. 
 \end{prop}
 \begin{proof} 
 By Proposition \ref{L3} we know that
$$ 	V^{1/3}\left( \int_{M}|(\frac{2}{3}s_{g}+ {2}w_{g} )_{-}|^{3}d\mu 
 \right)^{2/3}\geq 32\pi^{2} (c_{1}^{+})^{2}.$$
 On the other hand, 
 $$-4\inf \left(\underline{K}_{g} + \frac{s_{g}}{12}\right) \geq
 \left(\frac{2}{3}s_{g}+2w_{g}\right)_{-} $$
 at all points by virtue of $M$ by Lemma \ref{wgl}.
 Hence
 $$ V^{1/3}\left( \int_{M}\left|4 
 \inf \left(\underline{K}_{g} + \frac{s_{g}}{12}\right)\right|^{3}d\mu
 \right)^{2/3}\geq 32\pi^{2} (c_{1}^{+})^{2},$$
 and so
 $$
  \left|\inf \left(\underline{K}_{g} + \frac{s_{g}}{12}\right)\right|^{2}V
 \geq 2\pi^{2} (c_{1}^{+})^{2}.
 $$
 Taking the square root of this inequality then yields the claim.  
 \end{proof}

Let us now consider a  natural family of
 minimal volume invariants \cite{grom,bcg2,lno} defined with respect 
 to 
 only {\em lower bounds} on curvature. Let $M$ be any smooth compact 
 manifold of dimension $n$. By  just considering large constant
 multiples of  some given
 metric, one obtains metrics  on  on $M$ with  sectional curvature
 $K \geq -1$,  at the price of  making the total volume 
 of $M$ enormous. This  motivated Gromov \cite{grom} to define 
 the 
  {\em Gromov minimal volume}
 $$\mv_{K}(M):= \inf ~\left\{ ~\mv (M,g) ~|~ K_{g}\geq -1\right\}$$
as a way of quantifying the degree to which some negative sectional 
curvature might be an inevitable feature of all possible geometry
on certain manifolds. However, it is also very natural to
weaken the curvature hypothesis by instead only imposing 
lower bounds on the Ricci or  scalar curvatures. In particular, 
one may consider \cite{lno} what we might call the 
{\em Yamabe minimal volume} of $M$, defined by 
$$\mv_{s}(M):= \inf ~\left\{ ~\mv (M,g) ~\left| ~ \frac{s_{g}}{n(n-1)}
\geq -1\right. \right\},$$
where our conventions have been chosen so that,
by definition, 
$$\mv_{K}(M)\geq \mv_{s}(M)\geq 0.$$
One key advantage of the latter definition is
 that the mature theory of the Yamabe problem
 \cite{sch}
allows this minimal volume to be reinterpreted as 
$$\mv_{s}(M) = \inf_{g} 
\int_{M}\left|\frac{s_{-}}{n(n-1)}\right|^{n/2}d\mu_{g}
 = \left\{ \begin{array}{cl}0, & Y(M) \geq 0\\
\left(\frac{|Y(M)|}{n(n-1)}\right)^{n/2} , & Y(M) \leq 0,
\end{array}
\right.
$$
where $s_{-}=\min (s, 0)$, and where 
$Y(M)$ is the  Yamabe invariant (sigma constant) of $M$.

It will now be convenient to consider an interpolation
between these two definitions, gotten by averaging the 
two curvature conditions under consideration.
We thus define the {\em mixed minimal volume} of a 
the smooth compact manifold $M$ to be 
$$\mv_{K,s}(M) = 
\inf ~\left\{ ~\mv (M,g) ~\left| ~ 
\frac{1}{2}\left( K_{g}+ \frac{s_{g}}{n(n-1)}
\right)\geq -1\right. \right\}.$$
This definition has been chosen in such a way that we 
now have  the following result:

\begin{thm}\label{minvol}
Let $M$ be the underlying $4$-manifold of any complex
surface of general type. Then
$$\mv_{K,s}(M) \geq \frac{9}{4}\mv_{s}(M) .$$
Moreover, equality holds if  $M$ is any complex-hyperbolic manifold
${\mathbb   C}{\mathcal H}_{2}/\Gamma $, 
and in this case both minimal volumes are achieved by 
 appropriate constant multiples of the standard metric. 
\end{thm}
\begin{proof}
Let $X$ be the minimal model of $M$. Then for every metric $g$ on
$M$ there is a spin$^{c}$ structure 
with $(c_{1}^{+})^{2}\geq c_{1}^{2}(X)$ and such that
the corresponding Seiberg-Witten invariant is non-zero for
the polarization determined by $g$; in particular, 
the scalar curvature $s$ of $g$ is negative somewhere.
Thus,  Proposition \ref{lower} 
tells us that $g$ satisfies 
$$V \left[\min_{x\in M} \left( \underline{K}(x)+\frac{s(x)}{12}
\right) \right]^{2} 
 \geq 2\pi^{2}(c_{1}^{+})^{2}\geq  2\pi^{2}c_{1}^{2}(X) .
 $$
 Restricting our attention to metrics for which 
 $\frac{1}{2}(\underline{K}+ \frac{s}{12})\geq -1$, we thus get 
 $$\mv_{K,s}(M)\geq \frac{2\pi^{2}c_{1}^{2}(X)}{2^{2}}= 
 \frac{\pi^{2}}{2} c_{1}^{2}(X).$$
 However \cite{lno}, $\mv_{s}(M)= \frac{2\pi^{2}}{9} c_{1}^{2}(X),$
 so it follows that 
 $$\mv_{K,s}(M)\geq \frac{9}{4}\mv_{s}(M),$$
 as claimed. 
 
 If $M$ is a complex-hyperbolic manifold, 
 and if $g$ is a multiple of the  hyperbolic metric, normalized
 so that its sectional curvatures satisfy
 $K\in [-\frac{1}{3}, -\frac{4}{3}]$,  one then has
 $s=-8$, and $\frac{1}{2}(\underline{K}+s)\equiv -1$, whereas
 $s^{2}V= 32\pi^{2}c_{1}^{2}(M)= 32(\frac{9}{2} \mv_{s}(M))$. 
 Thus $V= \frac{9}{4} \mv_{s}(M)$, and 
 $$\mv_{K,s}(M)= \frac{9}{4}\mv_{s}(M)$$
 for any complex-hyperbolic $4$-manifold. 
\end{proof}

Because of the tautological inequality 
$$\mv_{K}(M)\geq 
\mv_{K,s}(M),$$ 
this immediately implies: 

\begin{cor}
If $M$ is any compact complex surface of general type, then 
$$\mv_{K}(M)\geq \frac{9}{4} \mv_{s}(M).$$ 
\end{cor}

Of course, the constant appearing in this corollary 
should not be expected to be sharp. Indeed, 
inspection of the complex-hyperbolic case 
would instead seem to 
suggest the following:

\begin{conj}
If $M$ is any compact complex surface of general type, then 
$$\mv_{K}(M)\geq 4 \mv_{s}(M),$$
with equality iff $M$  is complex hyperbolic.
\end{conj}

Notice that we have not even shown that 
only the complex hyperbolic manifolds saturate the
inequality of Theorem \ref{minvol}. Nonetheless, 
we can at least say the following:

\begin{prop}
Suppose that $M$ is a complex surface of general type, 
and suppose that the mixed 
minimal volume $\mv_{K,s}(M)$ is actually achieved 
by some metric $g$. Then either $\mv_{K,s}(M) > \frac{9}{4}\mv_{s}(M)$, 
or else $(M,g)$ is a complex-hyperbolic manifold,
 normalized so that 
$s\equiv -8$.
\end{prop}

\begin{proof}
Suppose that $g$  is a  metric with
${K}+\frac{s}{12}\geq	-2$ and  total volume 
$V=\frac{9}{4}\mv_{s}(M)$. Thus, letting 
$X$ denote the minimal model of $M$, 
and choosing a spin$^{c}$ structure on 
$M$ as in the proof of Proposition \ref{west}, 
\begin{eqnarray}
2\pi^{2}c_{1}^{2}(X)=
4 V	 & \geq & 
V\left[\min \left( \underline{K} +\frac{s}{12} \right)\right]^{2}
	\nonumber  \\
	 & \geq &  \frac{1}{16}\int (\frac{2}{3}s+ 2w)_{-}^{2}
 	d\mu 
	\label{aha}  \\
	 &  \geq & \frac{32 \pi^{2}}{16} [c_{1}^{+}]^{2}
	 = 2 \pi^{2}[c_{1}^{+}]^{2}
	\nonumber
\end{eqnarray}
 	by virtue of (\ref{yes}). Since the spin$^{c}$ structure 
 	 is chosen so that 
 	$[c_{1}^{+}]^{2}\geq c_{1}^{2}(X)$, equality must hold,
 	 it follows that $c_{1}^{+}= c_{1}(X)$ and  
 	$c_{1}^{-}=0$;  in particular, $M=X$ is minimal.  
 	 Moreover, since 
 	 (\ref{yes}) is  saturated,
 	 Proposition \ref{aka} tells us that $g$ is
 	almost-K\"ahler, with $s+s^{*}$ constant,  $\omega$  an  
  eigenvector of $W_{+}$ at each point,
 	and 
 	$c_{1}^{+}\propto [\omega ]$.
 	Since $c_{1}=c_{1}^{+}$, the almost-complex structure
 	satisfies $c_{1} \propto [\omega ]$, 
 	so that our almost-K\"ahler manifold is
 	{\em monotonic} in the terminology of \cite{td}.

Since equality holds in (\ref{aha}), we also have
$$\underline{K} +\frac{s}{12}\equiv \frac{s}{6}+ \frac{w}{2},$$
and this, for starters, tells us that $W_{-}\equiv 0$
by (\ref{crafty}). Inspection of (\ref{sec}) then 
shows that $B^{*}(\omega)\equiv 0$, so that 
$\omega$ is in fact an eigenvector of the 
full curvature operator ${\mathcal R}$. An almost-K\"ahler
manifold with this property is called {\em weakly $\ast$-Einstein}.
Because we also know that $(M,g,\omega )$ 
is monotonic, and that the sum  $s+s^{*}$ of
the scalar and $\ast$-scalar curvatures is constant, 
a result of Dr{\u{a}}ghici 
\cite[Proposition 2]{td}
asserts that $(M,g, J )$ is K\"ahler-Einstein. 
But we also have also observed that $W_{-}\equiv 0$,
so the entire curvature operator is parallel,
and the universal cover of $(M,g)$ is 
therefore the symmetric space ${\mathbb C}{\mathcal  H}_{2}$,
equipped with the unique multiple of its standard metric for which
$s=-8$.    
Thus $M$ admits a complex-hyperbolic metric; and this metric
is moreover unique, up to rescalings and diffeomorphisms, by
virtue of  Mostow rigidity. 
\end{proof}

We now close with some final remarks concerning minimal volumes. 
First of all, one might want to consider more general mixed minimal 
volumes, where the scalar and sectional curvatures are 
weighted differently from   above. It may therefore be
worth noting that one can also show that 
$$
\inf \left\{ \mv (M,g) ~\left| ~ tK_{g}+ (t-1) \frac{s}{12} \geq -1
\right. \right\} \geq (1+ t)^{2} \mv_{s} (M)
$$
for any complex algebraic   surface $M$ and any constant $t\in 
[0,\frac{1}{2}]$; and again,  equality  holds for  complex-hyperbolic
manifolds. Indeed, it is not hard to establish 
inequalities such as 
$$
\inf_{g}\int_{M}\left((1-\frac{2t}{3}) s + 4t w\right)_{-}^{2}d\mu
\geq 32\pi^{2} (c_{1}^{+})^{2}
$$
for $t\in [0,\frac{1}{2}]$, and the assertion then follows. The 
interesting question,
though, is whether such inequalities still hold for some
  $t> \frac{1}{2}$.

It would also be interesting to determine whether  
$$\mv_{r} (M) = \inf \left\{ \mv (M,g) ~\left| ~ r \geq - 3 g
\right. \right\} $$
is strictly larger than $\mv_{s} (M)$ when $M$ is a non-minimal 
complex surface; our computations of 
${\mathcal I}_{r}(M)$ certainly prove such an assertion for
the related invariant
$$
\inf \left\{ \mv (M,g) ~\left| ~3g \geq  r \geq - 3 g
\right. \right\}\geq \frac{1}{36}{\mathcal I}_{r}(M) ,
$$
defined with respect to {\em two-sided} curvature bounds. 
Because the minimizing sequences  constructed in the proof of 
Theorem \ref{l2ric} actually  have huge amounts of Ricci curvature concentrated
near the blow-ups,  both of these minimal volumes
will presumably turn out to be much larger than ${\mathcal I}_{r}/36$ 
in the non-minimal case. 
 A confirmation of  this speculation 
would presumably also lead to  
 even stronger non-existence theorems for Einstein metrics. 
And then, of course,  the tautological inequality
$$
\mv_{K} (M)\geq 
\mv_{r} (M)
$$
would most pointedly inform us that our present estimates of 
the  Gromov minimal volume have barely begun to scratch  the surface
of the subject.


\begin{thebibliography}{99}

\bibitem{ad} V. Apostolov, J.  Armstrong,  and T. Dr{\u{a}}ghici,
{\em Local Rigidity of Certain Classes of Almost K\"ahler 
$4$-Manifolds},
e-print math.DG/9911197, available at http://xxx.lanl.gov

\bibitem{rmstg}  J.  Armstrong, {\em An Ansatz for
Almost-K\"ahler, Einstein $4$-Manifolds},
Oxford preprint, 1999. 
       

\bibitem{aubin} T. Aubin, {\em
Equations du Type {M}onge-{A}mp\`{e}re sur les Vari\'{e}t\'{e}s 
{K\"a}hl{\'e}riennes  Compactes},
{\bf C. R. Acad. Sci. Paris 283A} (1976)   119--121.



\bibitem{aubin0} T. Aubin,  {\bf Nonlinear Analysis 
on Manifolds. Monge-Amp\`ere Equations},
Springer-Verlag, 1982. 


\bibitem{bpv} 
W. Barth,   C. Peters,  and A. Van de Ven,  
{\bf Compact Complex Surfaces},
Springer-Verlag, 
1984. 



\bibitem{bes} A. Besse, {\bf Einstein Manifolds}, Springer-Verlag,
1987.


\bibitem{bcg2} G. Besson, G. Courtois, and S. Gallot, 
{\em Entropies et Rigidit\'es des Espaces Localement Sym\'etriques de
Courbure Strictement N\'egative}, 
{\bf Geom. and Func. An. 5} (1995) 731--799.

\bibitem{blair} D.   Blair, 
{\em The `Total Scalar Curvature' as a Symplectic Invariant, and
             Related Results},
{\bf Proceedings of the 3rd Congress of Geometry 
(Thessaloniki,          1991)},
Aristotle Univ. Thessaloniki, 1992, pp. 79--83. 


\bibitem{bourg} J.P. Bourguignon, {\em  Les Vari{\'e}t{\'e}s de Dimension 
$4$ \`a Signature Non Nulle dont la Courbure
Est Harmonique Sont d'Einstein}, {\bf Invent.\  
Math.\  63}
 (1981) 263--286.
 
 \bibitem{chen} Z. Chen,
    {\em On the Geography of Surfaces. {S}imply Connected Minimal
             Surfaces with Positive Index},
  {\bf Math. Ann. 277} (1987) 41--164.
  
  
 
 \bibitem{rio} H. del Rio Guerra, {\em Seiberg-Witten Invariants of 
 Non-Simple Type and Einstein Metrics}, 
 e-print  math.DG/0002243, available at http://xxx.lanl.gov 
 
 
 \bibitem{td} T. Dr{\u{a}}ghici, 
   {\em On Some $4$-Dimensional Almost {K}\"ahler Manifolds},
  {\bf Kodai Math. J. 18} (1995) 156--168. 


\bibitem{freedman}
M. Freedman, {\em On the Topology of 4-Manifolds},
{\bf J. Diff. Geom. 17} (1982) 357--454. 
 
 
 
\bibitem{fs} R.  
Fintushel and R.  Stern, {\em Immersed Spheres in $4$-Manifolds and the 
Immersed
Thom Conjecture}, {\bf Turk.\  J. Math.\  19} (1995)  145--157.

\bibitem{grom} M. Gromov, {\em Volume and Bounded Cohomology},
{\bf Publ.\ IHES 56} (1982) 5--99.


\bibitem{gromsym} M. Gromov,
{\em Pseudoholomorphic Curves in Symplectic Manifolds},
{\bf Invent. Math.  82}  (1985) 307--347.

\bibitem{gvln} M.  Gromov and H. B.  Lawson,
{\em The Classification of Simply Connected Manifolds of Positive
             Scalar Curvature}, {\bf Ann. of Math. 111} (1980) 423--434.

\bibitem{gu} M. Gursky,  
{\em  Four-Manifolds with $\delta {W}^+=0$ and {E}instein  
Constants on the Sphere}, Indiana University preprint, 1997.

\bibitem{hit} N. Hitchin, {\em Compact Four-Dimensional Einstein 
Manifolds},
{\bf  J. Diff.\ Geom.\  9}
(1974) 435--441.



\bibitem{hori} E.  Horikawa, {\em Algebraic Surfaces of General Type 
with Small $c_{1}^{2}$,  I}, {\bf  Ann.\  Math.\  104} (1976) 357--387. 


\bibitem{rkob} R. Kobayashi, {\em Einstein-K\"ahler $V$-Metrics on 
Open Satake $V$-Surfaces with Isolated Quotient Singularities},
{\bf Math.\ Ann.\ 272} (1985) 385--398.

\bibitem{kot} D. Kotschick,
{\em  Einstein Metrics and Smooth Structures},
{\bf  Geom.\ Topol.\ 2} (1998) 1--10.




\bibitem{kron} P. Kronheimer, {\em Instantons Gravitationelles et
Singularit\'es de Klein}, {\bf C.R. Acad.\ Sci.\ Paris 303} 
(1986)  53--55.

 \bibitem{KM}  P. Kronheimer and T. Mrowka, 
{\em The Genus of Embedded Surfaces in the
Complex Projective Plane}, 
{\bf  Math.\ Res.\ Lett.\ 1}  (1994) 797--808.



 \bibitem{K}  P. Kronheimer,
 {\em Minimal Genus in $S^{1}\times M^{3}$}, 
{\bf Invent. Math. 135} (1999) 45--61. 
 

%% 
 % \bibitem{lrs}  C. LeBrun,
 % {\em	Scalar-Flat	K\"ahler	 Metrics	on Blown-Up	Ruled Surfaces}, 
 % {\bf	J. reine  angew.\ Math.\ 420} (1991) 161--177. 
 %%



\bibitem{lmo} C. LeBrun, {\em Einstein Metrics and Mostow Rigidity},
{\bf  Math.\ Res.\ Lett.\ 2} (1995)
1--8.
 

\bibitem{lpm} C. LeBrun, {\em Polarized 4-Manifolds, Extremal
K\"ahler Metrics, and
Seiberg-Witten Theory},  {\bf  Math.\ Res. Lett.\ 2} (1995)
653-662.



\bibitem{lno} C. LeBrun, 
{\em Four-Manifolds without Einstein Metrics},
{\bf Math.\ Res.\ Lett.
3} {(1996) 133--147.}


\bibitem{lky}  C. LeBrun, 
{\em Kodaira Dimension and the Yamabe Problem}, 
{\bf Comm. Anal. Geom. 7} (1999) 133--156.



\bibitem{lebweyl} C. LeBrun, {\em Weyl Curvature, 
{E}instein Metrics, and {S}eiberg-{W}itten
             Theory}, {\bf Math.\ Res.\ Lett.\  5} (1998)  423--438.


\bibitem{lebyond} C. LeBrun,
{\em 
Four-Dimensional 
Einstein Manifolds,
and Beyond},
to appear in {\bf 
Essays on Einstein Manifolds}, 
International Press, 2000. 




\bibitem{ozsz}
{P.~Ozsv\'ath and Z.~Szab\'o}, {\em Higher-Type Adjunction Inequalities in
  {S}eiberg-{W}itten Theory}.
\newblock Princeton University preprint, 1998.



\bibitem{jp}
{J.~Petean}, {\em Computations of the {Y}amabe Invariant}, {\bf Math. Res.
  Lett. 5} (1998) 703--709.



\bibitem{jp3}  J. Petean,  
{\em Yamabe Invariants of Simply Connected Manifolds}, 
e-print  math.DG/9808062, 
available at http://xxx.lanl.gov


\bibitem{samba} A. Sambusetti, {\em An Obstruction to the Existence of 
Einstein Metrics on 4-Manifolds}, {\bf  C. R.
Acad. Sci. Paris  322} (1996) 1213--1218. 


\bibitem{sch} R. Schoen, {\em  Variational Theory for the
Total Scalar Curvature Functional for {R}iemannian Metrics 
and Related Topics},
{\bf Lec. Notes Math. 1365} (1987)  120--154. 



\bibitem{st} 
I. M. Singer  and J. A.  Thorpe,  
{\em The Curvature of $4$-dimensional {E}instein Spaces},  
 in {\bf Global Analysis (Papers in Honor of K. {K}odaira)},
 D. Spencer, eds, 
 Univ. Tokyo Press,  1969, pp. 355--365. 
 

\bibitem{stolz} S. Stolz,
{\em Simply Connected Manifolds of Positive Scalar Curvature},
{\bf Ann. of Math. 136} (1992) 511--540.


 \bibitem{tasd}   C.H. Taubes, {\em The Existence of 
 Anti-Self-Dual Conformal Structures}, {\bf  J. Diff.
Geom. 36} (1992)
 163--253.
 
 \bibitem{taubes}   C.H. Taubes, {\em The Seiberg-Witten Invariants and 
Symplectic Forms}, {\bf Math.\ Res.\ Lett.\ 1} (1994) 809--822.



\bibitem{t2} C.H. Taubes, {\em The Seiberg-Witten and 
Gromov Invariants}, {\bf Math.\ Res.\ Lett.\ 2}
(1995) 221--238.

\bibitem{thorpe} J. Thorpe,
{\em Some Remarks on the Gauss-Bonnet Integral},
{\bf  J. Math.\  Mech.\  18} (1969) 779--786. 

\bibitem{trud} 
N. Trudinger,  {\em Remarks Concerning the
Conformal Deformation of Metrics to Constant Scalar Curvature}, 
{\bf Ann. Scuola Norm. Sup. Pisa 22} (1968) 265--274.

 
\bibitem{tsuj} H. Tsuji,
{\em Existence and Degeneration of K\"ahler-Einstein
Metrics on Minimal Algebraic Varieties of General Type},
{\bf Math. Ann. 281} (1988)  123--133. 

\bibitem{witten} E. Witten, {\em Monopoles and Four-Manifolds},
{\bf  Math.\ Res.\ Lett.\ 1}  (1994) 809--822.

\bibitem{yam} 
H. {Y}amabe,  
{\em On the Deformation of 
 {R}iemannian Structures on  Compact Manifolds},
 {\bf Osaka Math. J. 12} (1960) 21--37.
   

\bibitem{yau}  S.-T. Yau, {\em On the Ricci-Curvature of a Complex 
K\"ahler Manifold 
and the Complex Monge-Amp\`ere Equations}, {\bf Comment. Pure Appl.
Math.\ 31} 
(1978) 339--411.


\end{thebibliography}
\end{document}